\documentclass{amsart}
\usepackage{graphicx}
\usepackage{amssymb}
\usepackage{graphicx}
\usepackage{amsfonts}
\usepackage{amsmath}
\usepackage{amsthm}
\usepackage{esint}
\usepackage{mathabx}
\usepackage{mathtools}
\usepackage{nccmath}
\usepackage{enumitem}
\usepackage{romannum}
\usepackage[colorlinks]{hyperref}

\newtheorem{theorem}{Theorem}[section]
\newtheorem{lemma}{Lemma}[section]
\newtheorem{rmk}{Remark}[section]
\newtheorem{corollary}{Corollary}[section]

\newcommand{\bigO}{\mathcal{O}}
\newcommand{\e}{\epsilon}

\newcommand{\ic}{\int_{\mathbb{S}^1}}

\begin{document}
\pagenumbering{arabic}
\title[Regularity for the equation $\det(D^2u) = f(u)$]{Regularity theory for the equation $\det(D^2u) = f(u)$ in two dimensions}

\author{Connor Mooney}
\address{Department of Mathematics, UC Irvine}
\email{\tt mooneycr@uci.edu}

\author{Aranya Sen}
\address{Department of Mathematics, UC Irvine}
\email{\tt aranyas1@uci.edu}


\begin{abstract}
We study the interior regularity of convex solutions to
$$\det D^2u = f(u)$$
in planar domains, with no regularity or monotonicity assumptions on $f$. In particular, when $f$ is bounded between positive constants, we show that $u$ is locally $C^{1,\,1}$, and that its level sets near a critical point are asymptotically ellipses. Our analysis relies on a new monotone quantity for a second-order flow satisfied by the support functions of the level sets of $u$. As a byproduct of our methods, we obtain Liouville theorems that generalize a classical result of J\"{o}rgens, and include a complete classification of global solutions in 2D to an obstacle problem studied by Jin-Tu-Xiong. As another application, we obtain regularity and rigidity results for rotating-wave solutions of the semigeostrophic equations from meteorology.
\end{abstract}
\maketitle


\section{Introduction}
In this paper we study the interior regularity of convex solutions to
\begin{equation}\label{MainMA}
\det D^2u = f(u)
\end{equation}
in the plane, with no regularity or monotonicity assumptions on $f$. Our interest in (\ref{MainMA}) stems from the fact that the solutions characterize the rotating-wave solutions (i.e. steady states, up to a constant-speed background rotation of the plane) of the semigeostrophic equations from meteorology. Equations of this type also arise in geometric contexts, such as the $L_p$ Minkowski problem and related free boundary problems (see \cite{JTX}, \cite{CW}, \cite{Lut}, \cite{LYZ}). Our results have interesting consequences for both applications.

We first state a simplified version of our main result. If 
$$\det D^2u \in [\lambda,\,\Lambda], \quad 0 < \lambda \leq \Lambda < \infty,$$ 
then Caffarelli's regularity theory implies that $u$ is locally $C^{1,\,\alpha}$ and strictly convex \cite{C1}. Work of De Philippis-Figalli \cite{DF}, De Philippis-Figalli-Savin \cite{DFS}, and Schmidt \cite{Schimdt} shows that $u$ is also locally $W^{2,\,1+\epsilon}$. This local regularity is optimal without further information \cite{Wang}. For equations with the special structure (\ref{MainMA}), our main result improves the regularity substantially:

\begin{theorem}\label{Main}
Assume that $f \in [\lambda,\,\Lambda]$, and that $u$ solves (\ref{MainMA}) in a normalized section $S_1(0) \subset \mathbb{R}^2$. Then
\begin{enumerate}[label=\Roman*.]
\item $u \in C^{1,\,1}_{loc}$ and
\begin{equation}\label{C11}
\|D^2u\|_{L^{\infty}(S_{1/2}(0))} \leq C(\lambda,\,\Lambda).
\end{equation}
\item For $t > 0$, the level sets $\{u = t\}$ are smooth and have precise asymptotic expansions in terms of $f$ as $t$ varies. 
\item As $t \rightarrow 0$, the area-normalized level sets $\frac{\{u = t\}}{|\{u < t\}|^{1/2}}$ converge smoothly to an ellipse.
\end{enumerate}
\end{theorem}
We can in fact say the same under weaker assumptions that allow $f$ to degenerate in a controlled way (see Remark \ref{rmk:gen}), and we can prove both \Romannum{2}. and \Romannum{3}. while still allowing $f$ to blow up in a controlled way (see Remark \ref{rmk:gen} and Theorem \ref{x=0_theorem}). We also believe that the level sets at positive height are analytic, and we leave this for future work.

Important features of our result are that no assumptions about the regularity or monotonicity of $f$, or about the regularity of $\partial S_1(0)$, are made. There has been a great deal of interest in bounding the curvature of the level sets of solutions to (\ref{MainMA}), see e.g. \cite{CMS}, \cite{FS}, \cite{HHW}, \cite{SL}. In \cite{CMS} and \cite{FS} the authors identify quantities that satisfy the maximum principle when $f = 1$ and $f$ is monotone increasing, respectively, which are only effective when assumptions about boundary curvatures are made. In \cite{HHW} and \cite{SL} the authors study the Monge-Amp\`{e}re eigenvalue problem in smooth and uniformly convex domains in 2D and higher dimensions, respectively, and the primary concern is boundary regularity (and in particular, getting lower bounds on level set curvatures near the boundary). 

In contrast with the proofs of the above-mentioned results, which are mainly maximum principle-based, the proof of Theorem \ref{Main} relies on a new monotone quantity for a second-order flow that is satisfied by the support functions of the level sets (discussed in greater detail below). As a byproduct of our methods, we obtain several Liouville theorems. For example, we show that global solutions to (\ref{MainMA}) must be radial, up to an affine transformation, under relatively mild conditions on $f$:
\begin{theorem}\label{Liouville2}
Let $u$ be a nonnegative solution to (\ref{MainMA}) on $\mathbb{R}^2$ with an isolated minimum. Assume further that
$$
\lambda t^p\le f(t)\le \Lambda t^p
$$
for some $p \in (-2,\,2)$. Then up to an affine change of variable, $u$ is radial.
\end{theorem}
\noindent When $f = 1$, Theorem \ref{Liouville2} is classical and is due to J\"{o}rgens \cite{Jorg}  (and was extended to dimension five by Calabi \cite{Cal} and to all dimensions by Pogorelov \cite{Pog}).

We can also classify the global nonnegative solutions on $\mathbb{R}^2$ to the free boundary problem
\begin{equation}\label{FB}
\det D^2u = u^p\chi_{\{u > 0\}}, \quad 0 \leq p < 2
\end{equation}
with bounded coincidence set:
\begin{theorem}\label{Main2}
Global solutions to (\ref{FB}) in $\mathbb{R}^2$ with bounded coincidence set $\{u  = 0\}$ are radial, up to an affine transformation.
\end{theorem}
\noindent Theorem \ref{Main2} completes in 2D a program initiated by Jin-Tu-Xiong \cite{JTX}, which was motivated by applications to the $L_p$ Minkowski problem.  

Finally, we can classify the global solutions to certain degenerate Monge-Amp\`{e}re equations in the plane that depend only on one variable:
\begin{theorem}\label{Liouville1}
Let $u$ be a solution on $\mathbb{R}^2$ to
$$\det D^2u = g(x_2),$$
where $\lambda |t|^m \leq g(t) \leq \Lambda |t|^m$ for some $m > -1$. Then up to affine changes,
$$u = \frac{a}{2}x_1^2 + \frac{1}{a}\int_0^{x_2}g(t)(x_2-t)\,dt+l(x_1,x_2)$$
for some constant $a>0$ and a linear function $l$.
\end{theorem}
\noindent Theorem \ref{Liouville1} generalizes the main result in \cite{JX}, which deals with the case $g(t) = |t|^m$.

Other interesting consequences of our main result concern the semigeostrophic equations (SGEs), which model large-scale atmospheric and oceanic flow. The SGEs, in so-called dual coordinates, are a system for a spatially convex function $u$ and a density $\rho \geq 0$ on $\mathbb{R}^2 \times [0,\,T)$. After a constant-speed rotation of the plane, the dual SGEs can be written as
\begin{equation}\label{DSGEs}
\begin{cases}
\partial_t\rho - \text{div}(\rho \nabla^{\perp}u) = 0, \\
\det D^2u = \rho.
\end{cases}
\end{equation}
In particular, the dual SGEs are the same as 2D Euler, with the Monge-Amp\`{e}re equation replacing the Laplace equation. Because the level sets of $u$ are connected whenever $\rho$ is nonzero, the steady states of (\ref{DSGEs}) are exactly the solutions to (\ref{MainMA}), with $f = \rho$.  This distinguishes (\ref{DSGEs}) from 2D Euler, which admits steady states that don't arise from solutions to $\Delta u = \rho(u)$ unless one imposes additional hypotheses (e.g. analyticity, see \cite{EHSX}, \cite{EFR}, and \cite{GSPS}).
More remarkably, our main result implies that for steady state with $\rho$ locally bounded between positive constants, $u \in C^{1,\,1}_{loc}$. This also stands in contrast with 2D Euler, which admits e.g. vortex patch steady states with unbounded Hessian such as the Bahouri-Chemin example in \cite{BC}. As a consequence of this regularity result, we are able to show that the physical velocity of rotating-wave solutions to the SGEs is locally bounded (see Theorem \ref{Thm_SGE1}). We are also able to establish rigidity results for global rotating-wave solutions, under relatively mild assumptions on $\rho$, which say that the trajectories are ellipses, after a constant speed rotation of the plane. See Theorem \ref{SGE_thm2} for a precise statement, and Section \ref{semigeostrophic} for figures. Since steady states often capture long-time behavior of solutions to dynamical systems, one might view Theorem \ref{Main} as evidence that the long-time behavior of solutions to the SGEs is in fact quite regular.

\begin{rmk}
This corroborates findings from \cite{Cullen} which suggest that the SGEs may be less turbulent than Euler. This observation was also noted in \cite{MO}, where the authors construct interesting solutions to the SGEs with pressure and stream function quadratic at each fixed time.
\end{rmk}

The majority of the work in this paper goes into proving Theorem \ref{Main}. Its proof consists of two main parts. The first part addresses \Romannum{2}., namely, the regularity and expansions for $u$ near its nondegenerate level sets. For this we use the Hodograph transform, the partial Legendre transform, and Schauder-type estimates. The argument benefits greatly from the framework developed by Daskalopoulos-Savin in \cite{SD} for Monge-Amp\`{e}re equations with point degeneracies. It is similar in spirit to the work in \cite{Li} for Monge-Amp\`{e}re equations with right hand sides that are only regular in one variable, but gives a more precise pointwise Schauder-type expansion. The second part addresses \Romannum{3}. (from which \Romannum{1}. quickly follows). For this we study a second-order flow for the support functions of the level sets of $u$. The key is the discovery of an affine-invariant monotone quantity whose dissipation rate controls how close the level sets are to evolving by homothety and to being ellipses. Our analysis also relies on the variational structure of this flow, whose affine invariance gives a useful conserved quantity. The proof of monotonicity combines the flow with classical mixed-volume inequalities. In the context of geometric flows, monotone quantities have been derived in a similar way (see e.g. Andrews in \cite{A1} and \cite{A3}, Sapiro-Tannenbaum \cite{ST}).

The paper is organized as follows. In Section \ref{sec:Pre} we recall some preliminary results about sections of solutions to the Monge-Amp\`{e}re equation and about the Hodograph and partial Legendre transforms. We also record a few well-known facts from convex analysis. In Section \ref{sec:expansion} we address \Romannum{2}., namely, the regularity near nondegenerate level sets and as a consequence, the $C^{1,\,1}$ regularity of $u$ away from $0$. In Section \ref{sec:Liouville} we use these results to prove Theorem \ref{Liouville1}. In Section \ref{sec:min_point_analysis} we address \Romannum{3}., that is, the asymptotic roundness (up to affine transformations) of level sets at low height, which combined with the prior results gives \Romannum{1}. A complete proof of Theorem \ref{Main} is also provided. In Section \ref{application} we discuss applications of this asymptotic roundness, in particular, Liouville theorems Theorem \ref{Liouville2} and Theorem \ref{Main2} as well as certain rigidity results for \eqref{MainMA}. In Section \ref{semigeostrophic} we discuss applications of Theorem \ref{Main} to rotating-wave solutions of the semigeostrophic equations. Finally, in Section \ref{OpenQuestions} we state some open questions for future investigation motivated by this work.

\section*{Acknowledgments}
C. Mooney was supported by NSF Grant DMS-2143668. A. Sen was also supported by the same grant. A. Sen would like to thank Nam Q. Le for suggesting the free boundary problem considered in Theorem \ref{Main2} and Chong Gu for valuable discussions.

\section*{AI Disclosure}
The monotone quantity for the flow studied in Section \ref{sec:min_point_analysis} was discovered with the help of ChatGPT Plus. The remaining work (including all writing) was done solely by the authors.



\section{Preliminaries}\label{sec:Pre}
\subsection{Doubling measure RHS to Monge-Ampere equations}

For a differentiable convex function $u: \Omega \to \mathbb{R}$, we denote the section of height $t$ at a point $x\in \Omega$ to be
\begin{equation}\label{section}
    S_t(u,x) = \{y\in \Omega \;|\; 0 \leq u(y)- u(x)-\nabla u(x)\cdot(y-x) \leq  t\}.
\end{equation}

We also define a doubling measure $\mu$ w.r.t. ellipsoids in $\Omega$ if there exists a constant $c>0$ such that for any point $x_0\in \Omega$ and an ellipse $E$ such that $x_0+E\subseteq \Omega$, we have,
$$\mu(x_0+E) \geq c\mu((x_0+2E)\cap \Omega).$$

 Some important results while dealing with doubling measures we use are due to Caffarelli in \cite{C1}, Caffarelli-Gutierrez in \cite{CG1}, Jerison in \cite{J1} and Gutierrez-Hartenstine in \cite{GH1}.
\begin{lemma}\label{balacing}
    Assume $\mu$ to be a doubling measure. Let $u:\Omega\to \mathbb{R}$ be an Alexandrov solution of
    $$\det(D^2u) = \mu.$$
    Then for a compactly contained section $S_h(u,x)$ there exists a matrix $A_h$ such that
    $$k_0^{-1}A_hB_{r_h} \subseteq S_h(u,x)-x\subseteq k_0A_hB_{r_h},$$
    for some constant $k_0$ depending only on the doubling constant of $\mu$ and universal constants. Moreover $r_h = h(\mu(S_h(u,x)))^{-\frac{1}{n}}$ and $\det(A_h) = 1$.
\end{lemma}

\begin{lemma}\label{SCR} Let $u:\Omega\to\mathbb{R}$ be a solution of
    $$\det(D^2u) = \mu,$$
    where $\mu$ is a doubling measure whose density w.r.t. Lebesgue measure is bounded. Then for a compactly contained section $S_1(u,x)$, for $y\in S_{1}(u,x)$ the inequalities
    $$L_x(y)+c_1\omega(|y-x|) \leq u(y) \leq L_x(y)+C_1|x-y|^{1+\alpha},$$
    hold for universal constants $c_1,C_1,\alpha$, a modulus of strict convexity, $\omega$, and $L_x$, the supporting hyperplane to $u$ at $x$.
\end{lemma}

\begin{lemma}\label{line_sing}
Assume $\mu$ to be a doubling measure. Let $u:\Omega\to \mathbb{R}$ be an Alexandrov solution of
    $$\det(D^2u) = \mu.$$
    Then for a linear function $L$, the set $\{u=L\}$ is either singleton or has no extreme points in $\Omega$.
\end{lemma}

\subsection{Convex analysis} Here we recall some basic convex analysis facts in $2$ dimensions (see for e.g. Schneider \cite{Sch}). The classic Brunn-Minkowski theory states that the square root of the volume is a concave functional on the set of convex bodies, i.e., for two convex bodies $K,L\subseteq \mathbb{R}^2$ and $0\leq s\leq 1$ we have,
$$|(1-s)K+sL|^{\frac{1}{2}} \geq (1-s)|K|^{\frac{1}{2}}+s|L|^{\frac{1}{2}}.$$
Moreover there is an exact expression in terms of mixed volume, $V(K,L)$,
$$|K+L| = |K|+|L|+2V(K,L).$$
The following Minkowski inequality follows,
$$V(K,L) \geq \sqrt{K}\sqrt{L}.$$
An important fact is that if $K,L$ were two convex bodies with non empty interior then there is equality in the above Minkowski inequality if and only if $K$ and $L$ are homothetic to one another.
We have the following lemma for convex functions.
\begin{lemma}
    Let $u:S_1(u,0) \to \mathbb{R}$ be a convex function. Then for $A(t) = |S_t(u,0)|$, $\sqrt{A(t)}$ is a concave function.
\end{lemma}
\begin{proof}
    Assume without loss of generality $u(0) = 0$ and $\nabla u(0) = 0$. We note that since $u$ is convex, for $0\leq s\leq 1$ and $0\leq t_1,t_2\leq 1$,
    $$(1-s)S_{t_1}(u,0)+sS_{t_2}(u,0) \subseteq S_{(1-s)t_1+st_2}(u,0).$$
    Therefore by Brunn Minkowski inequality we have,
    $$(1-s)|S_{t_1}(u,0)|^{1/2}+s|S_{t_2}(u,0)|^{1/2} \leq |S_{(1-s)t_1+st_2}(u,0)|^{1/2}.$$
    This shows concavity of $\sqrt{A(t)}$.
\end{proof}
We also prove a basic growth lemma.
\begin{lemma}
    Let $f:\mathbb{R}_{\geq 0}\to \mathbb{R}_{\geq 0}$ be a differentiable function such that for some $0<  \alpha <2$,
    $$c t^{\alpha} \leq f(t) \leq Ct^{\alpha},$$
    If $\sqrt{f}$ is concave then, 
    $$\tilde{c}t^{\alpha-1}\leq f'(t) \leq \tilde{C}t^{\alpha-1}$$
    for some constants $\tilde{c},\tilde{C}$ depending only upon $c,C,\alpha$.
\end{lemma}
\begin{proof}
    Let $g = \sqrt{f}$ be concave. We have,
    $$\sqrt{c}t^{\alpha/2}\leq g \leq \sqrt{C}t^{\alpha/2}.$$
    By monotonicity of its derivative we have for $K$ large,
    $$\frac{g(Kt)-g(t)}{t(K-1)} \leq g'(t) \leq \frac{g(t) -g(0)}{t} .$$
    Using the uniform bounds on $g$ we get,
    $$\frac{\sqrt{c}K^{\alpha/2}-\sqrt{C}}{K-1}t^{\alpha/2-1} \leq g'(t) \leq \sqrt{C}t^{\alpha/2-1}.$$
    Therefore for $K$ large we get $g'(t)\sim t^{\alpha/2-1}$ (see Remark \ref{rmk:theta}). Now note that since $f = g^2$, we have $f' = 2g g'$. Therefore combining the estimates for $g$ and $g'$, the conclusion on $f'$ follows.
\end{proof}

\begin{rmk}\label{rmk:theta}
    Note that we call two positive functions $f\sim g$ if there exist constants $c,C$ such that,
    $$cg \leq f \leq Cg.$$
\end{rmk}

\subsection{Hodograph transform}\label{Hodograph}

We call $v$ the Hodograph transform (in the $y-$variable) of $u$ if it satisfies the following implicit relation
\begin{equation}\label{HT}
    u(x,v(x,z)) = z.
\end{equation}
Now assume $u: S_1(u,x)\to \mathbb{R}$ solves
\begin{equation}
    \det(D^2u) = \rho(u).
\end{equation}
By differentiating equation \eqref{HT} we get,
$$v_x = -\frac{u_x}{u_y},\;\;\text{ and }\; v_z = \frac{1}{u_y}.$$
We also note that the sets,
$$G_u = \{(x,y,z)\;|\; u(x,y) = z \}\;\text{ and }\; G_v = \{(x,y,z)\;|\;v(x,y) =z \}$$
are a reflection of one another about the plane $z = y$. Therefore by comparing Gaussian curvatures of both surfaces,
$$K_u = K_v \implies \det(D^2u) = \frac{\det(D^2v)}{v_y^4}.$$
Thus the hodograph transform satisfies the equation,
\begin{equation}\label{HTE}
    \det(D^2v) = v_z^4\rho(z).
\end{equation}
Note either $v$ or $-v$ will be convex and we choose appropriately. If $S_h(u,y) \subseteq S_1(u,x)\setminus S_\e(u,x)$ the hodograph transform $v$ would satisfy \eqref{HTE} in a section of height $\sim \sigma_\e h$ for some $\sigma_\e>0$. Therefore strict convexity at the minimum point ensures we can find universal $\sigma_\epsilon$ for each $\epsilon$ (See Theorem \ref{x=0_theorem}).

\subsection{Partial Legendre Transform}\label{subsec:PLT}

In this section we will study a model equation for a degenerate but doubling right hand side given by
\begin{equation}\label{doubling_equation}
    \det(D^2v)(x,z) = \rho_m(z).
\end{equation}
Here, for some $m> -1$, and $0<c_0\leq C_0$,
\begin{equation}\label{dr}
   c_0|z|^m\leq \rho_m(z)\leq C_0|z|^m.
\end{equation}
A model solution to a general right hand side $\rho$ can be given by
$$w_\rho(x,z) = \frac{|x|^2}{2}+\int_0^z\rho(t)(z-t)dt.$$
We define the following model solutions to equation \eqref{doubling_equation} as,
\begin{equation}
    w_m(x,z) = \frac{|x|^2}{2}+\int_0^z\rho_m(t)(z-t)dt.
\end{equation}
To reflect the natural scaling of the problem we also introduce the following anisotropic quasi-norm
 
\begin{equation}
    |(x,z)|_m = \sqrt{|x|^2+|z|^{m+2}}.
\end{equation}
We will now use the partial Legendre transform to convert the above Monge-Ampere equation into a linear partial differential equation given by,
\begin{equation}
v^{*}_{zz} + \rho_m(z)v^{*}_{xx} = 0.
\end{equation} 

We recall some important properties of partial Legendre transform discussed in \cite[Remark 2.3]{SD},
\begin{enumerate}
    \item Partial Legendre transform is involutive, i.e., $(u^{*})^{*}=u$.
    \item If $|u-v|\leq \e$ then $|u^*-v^*|\leq \e$ on their common domain of definition.
    \item For $e>0$ the partial Legendre transform of,
         $$g(x,z) = e\int_0^z\rho_m(t)(z-t)dt + \frac{(x-dz)^2}{2e}$$
         by keeping the $z-$variable constant is given by,
         $$g^*(p,z) = dpz+e\left(\frac{p^2}{2}-\int_0^z\rho_m(t)(z-t)dt\right).$$
\end{enumerate}

We consider the following $2^{nd}$ order expansion lemma,
\begin{lemma}\label{PLT_expansion}
    For $m\geq 0$ let $u$ be a solution to the following linear PDE in $B_1(0)$
    \begin{equation}\label{plt_lemma_eqn}
        L(u):= u_{zz} + \rho_m(z)u_{xx} = 0
    \end{equation}
    such that $\|u\|_{L^{\infty}(B_1(0))}\leq 1$. Then for $(x,z)\in B_{1/2}$ the following expansion holds
    $$u(x,z) = a+bx+cz+dxz+e\left(\frac{x^2}{2} - \int_0^z\rho_m(t)(z-t)dt\right)+\bigO(|(x,z)|_m^{2+\delta})$$
    for some $\delta>0$. Moreover the coefficients are bounded by universal constants.
\end{lemma}
\begin{proof}
    Note that $u_x,u_{xx},u_{xxx}$ also satisfy the equation \eqref{plt_lemma_eqn}. We now use the Bernstein technique to show boundedness of derivatives along $x$ direction. For a cutoff function $\phi$ (to be chosen later), consider
    $$w = Cu^2+\phi^2u_x^2.$$
    By a direct computation,
    \begin{equation}
        \begin{split}
            Lw = 2C(u_z^2+\rho_m(z)u_x^2)+L(\phi^2)u_x^2+ & (2\phi^2u_{xz}^2+8\phi u_{xz}\phi_zu_x)+\\
    & \rho_m(z)(2\phi^2u^2_{xx}+8\phi u_{xx}\phi_xu_x).
        \end{split}  
\end{equation} 
Hence,
\begin{equation}
    \begin{split}
        Lw \geq 2C\rho_m(z)u_x^2+L(\phi^2)u_x^2-32\phi_z^2u_x^2-32\rho_m(z)\phi_x^2u_x^2.
    \end{split}
\end{equation}
Now the cutoff function is chosen to be $\phi = F(x)F(z)$. Here $F$ is a one dimensional cutoff function such that $F([-1/4,1/4]) = 1$, $F([-1/2,1/2]^c)=0$ and $0\leq F\leq 1$. Therefore on $|z|\leq 1/4$, $\phi_z = 0$, and for a large $C$, $L(\phi^2)-32\phi_z^2\geq -C\rho_m(z)$. For this choice of $\phi$ and $C$
$$L(w)\geq 0.$$
 Thus by the maximum principle $u_{x}$ is bounded. By a similar computation $u_{xx},u_{xxx},u_{xxxx}$ are also bounded. Therefore by using the equation \eqref{plt_lemma_eqn} we can bound $u_{zz},u_{xzz},u_{xxzz}$. Now using the fact that for a single variable function $f$, if $f,f"$ are bounded, one can bound $f'$, we also bound $u_{z},u_{xz},u_{xxz}$. Now that $D^2u$ bounds have already been established, we analyze the remainder term in the Taylor expansion of $u$,
\begin{equation}\label{taylor_exp}
    \begin{split}
        u(x,z) &= u(x,0)+u_z(x,0)z+\int_0^z (z-t)u_{zz}(x,t)dt\\
        &= u(x,0)+u_z(x,0)z-u_{xx}(x,0)\int_0^z (z-t)\rho_m(t)dt+\bigO(|z|^{m+3}).
    \end{split}
\end{equation}
Further the Taylor expansions of $u(x,0),u_z(x,0)$ and $u_{xx}(x,0)$ are given by
$$u(x,0) = u(0,0)+u_x(0,0)x+u_{xx}(0,0)\frac{x^2}{2}+\bigO(|x|^3),$$
$$u_{z}(x,0) = u_z(0,0)+u_{xz}(0,0)x+\bigO(|x|^2),$$
$$u_{xx}(x,0) = u_{xx}(0,0)+\bigO(|x|).$$
Plugging these back in \eqref{taylor_exp},
 \begin{equation}
    \begin{split}
        &u(x,z) = u(0,0)+u_x(0,0)x+u_z(0,0)z+u_{xz}(0,0)xz+u_{xx}(0,0)\frac{x^2}{2}\\
        &-u_{xx}(0,0)\int (z-t)\rho_m(t)dt+\bigO(|z|^{m+3}+|x|^3+|x||z|^{m+2}+|x^2z|)
    \end{split}
\end{equation}
The mixed remainder terms can be bounded by Young's inequality.
\end{proof}

\begin{rmk}\label{rmk:neg_m}
    Note that for the case $0>m>-1$, the statement still holds if $u\in C^1$ across the line $z=0$. The only difference is that it is not immediate that $u_z,u_{xz},..$ are bounded. However we can bound $u_{zz},u_{xzz},...$ slightly away from the $z=0$ axis, say at $z=1/8$ and then use \eqref{plt_lemma_eqn} to integrate along lines parallel to x-axis. Moreover for $-1<m<0$, $\delta$ can be arranged in the statement of the lemma so that the error term is always superquadratic wrt $z$.
\end{rmk}

\section{Second order expansion and $C^2$ estimates}\label{sec:expansion}
In this section we study Monge-Ampere equations of the form,
\begin{equation}
    \det(D^2v)(x,z) = cf(x,z)\rho_m(z)\;\;,\;(x,z)\in S_1 
\end{equation}
 for some $f\in C^{\alpha}(S_1)$, such that $0<c_0\leq f\leq C_0$. We further assume $v\geq 0$ and $v(0,0) = 0$ and $\nabla v(0,0) = 0$. For such a solution we simplify the notation of a section to just $S_t = S(v,(0,0))$. For technical reasons we carry a constant $c$ on the right hand side with no apriori information about it. We call a constant universal if it depends on $m,c_0,C_0,\alpha$. Equations of the above form were studied by Li in \cite{Li} for $m=0$. Here we use a different proof strategy first used by Savin-Daskalopoulos in \cite{SD}.

\begin{lemma}
    \label{BC}
        Let $v$ in a section $S_1$ satisfy
    $$\det(D^2v)(x,z) = cf(x,z)\rho_m(z),$$
    $$S_{\gamma}(w_m,(0,0)) \subseteq S_1 \subseteq S_{\frac{1}{\gamma}}(w_m,(0,0)),$$
    for some constant $c>0$ and function $f(x,z)$. Then given $\theta_0$ small, there exist $\epsilon_0,t_0 > 0$ small such that if
    $$|f-f_0|\leq \epsilon_0$$
    for some constant $f_0>0$, then,
    $$A_0S_{t_0(1-\theta_0)}(w_{\rho_m\circ A_0},(0,0)) \subseteq S_{t_0} \subseteq A_0S_{t_0(1+\theta_0)}(w_{\rho_m\circ A_0},(0,0)),$$
    $$A_0 =  \begin{bmatrix}
        a_{11} & a_{12} \\
        0 & a_{22} 
    \end{bmatrix}, \;\;\; c_A\leq|a_{ii}|\leq C_A,\;i=1,2 \;\text{ and }\; |a_{12}|\leq C_A.$$
    Here $\epsilon_0,t_0,c_A,C_A$ only depend on $\gamma,\theta_0,f_0$ and universal constants.
    
\end{lemma}

\begin{proof}
    
  The constant $c$ is controlled by the geometry of section $S_1$, by comparing $\sqrt{cf_0}w_m$ and $v$.
    This can be done by considering the following inequalities,
    \begin{equation}
        \begin{split}
            \det (D^2(v+\sqrt{c\e_0}w_m))\geq \det&(D^2v)+\det(D^2(\sqrt{c\e_0}w_m))\\
            &\geq c(f+\e_0)\rho_m(z)\geq \det(D^2(\sqrt{cf_0}w_m)).
        \end{split}
    \end{equation}
    
    Now if on $\partial S_1$, $v+\sqrt{c\e_0}w_m < \sqrt{cf_0}w_m$. This leads to a contradiction of the maximum principle since $v(0,0) = w_m(0,0) = 0$. Thus an upper bound is established for $cf_0$ for small $\e_0$. Similar arguments with $\sqrt{cf_0}w_m+\sqrt{c\e_0}w_m$ and $v$, gives lower bound. Therefore,
    $$\frac{1}{C} \leq cf_0\leq C$$
    \\
    Let $h$ be a function such that
    $$\det(D^2h)(x,z) = cf_0\rho_m(z),\;\; h=1 \text{ on }\partial S_1.$$
    By a Pogorelov type estimate in the $x$ variable,
    $$0<\partial_{xx}h\leq C.$$
    Since $|f-f_0|\leq \epsilon_0$, by comparing $h+\sqrt{c\e_0}(w_m-\frac{1}{\delta})$ with $v$ and also $h$ with $v+\sqrt{c\e_0}(w_m-\frac{1}{\delta})$, the comparison principle implies,
    $$|h-v|\leq C\sqrt{\epsilon_0}\;\;\text{ and }\; |h|\leq 1.$$
     Taking partial Legendre transform of $h$ by fixing the $z$-variable, it satisfies
    $$\partial_{zz}h^*+ cf_0\rho_m(z)\partial_{pp}h^* = 0,$$
    and, 
    $$|h^*-v^*|\leq C\sqrt{\epsilon_0}\;\;\text{ and }\;|h^*|\leq 1.$$
    Thus by Lemma \ref{PLT_expansion},
    $$h^*(p,z) = a + bp+cz+dpz+e\left(\frac{p^2}{2}-cf_0\int_0^z\rho_m(t)(z-t)dt\right)+\bigO (|(p,z)|_m^{2+\delta}).$$

    Since $\partial_{xx}h \leq C$, $e = \partial_{pp}h^*(0,0)\geq \frac{1}{C}$, which guarantees uniform convexity in the $p-$variable.
    Also, $v(0,0) = 0,\nabla v(0,0)=0$. Therefore $v^*(0,0)=0, \nabla v^*(0,0) = 0$. Since $v^*$ is convex in the first variable,
    $v^*(p,0)\geq 0$, and by definition of partial Legendre transform, $v^*(0,z)\leq 0$. Thus for all $(p,0)$ and $(0,z)$ in the domain,
    $$h^*(p,0)\geq -C\sqrt{\epsilon_0},\;\;h^*(0,z)\leq C\sqrt{\epsilon_0}.$$
    This bounds the linear part of the expansion,
    $$|a| \leq C\e_0^{\frac{1}{2}},\;|b|,|c|\leq C\e_0^{\delta_m}$$
     for $\delta_m \leq\min\{\frac{1}{2},\frac{m+1}{2(m+2)}\}$. Thus for $|(p,z)|_m^2\leq \beta t_0$ (for some large $\beta$ depending on the uniform convexity of $h^*$ in the $p$ variable and smallness of $t_0$), we have
    \begin{equation}\label{lem1_exp}
    h^*(p,z) = dpz+e\left(\frac{p^2}{2}-cf_0\int_0^z\rho_m(t)(z-t)dt\right)+\bigO (t_0^{1+\frac{\delta}{2}}+\e_0^{\delta_m}).
    \end{equation}
    Define $g$ by,
     $$g(x,z) = e\int_0^zcf_0\rho_m(t)(z-t)dt + \frac{(x-dz)^2}{2e}.$$
    Then by equation \eqref{lem1_exp},
    $$|h^*-g^*|\leq C(t_0^{1+\frac{\delta}{2}}+\e_0^{\delta_m})$$
    Now since $g^*$ is uniformly convex in the $p$ variable, for $t_0$ small enough,
    \begin{equation}\label{error_bound_lem1}
    |h-g|\leq C(t_0^{1+\frac{\delta}{2}}+\e_0^{\delta_m}),
    \end{equation}
    given $|(\partial_xg,z)|_m^2 \leq \beta t_0$, i.e.,
    $$S = \{(x,z)\;|\;|z|^{m+2}+\frac{(x-dz)^2}{e^2}\leq \beta t_0\}$$
    Now let
    $$A_0 = \begin{bmatrix}
        \sqrt{e} && \frac{d}{\sqrt{cf_0e}}\\
        0 && \frac{1}{\sqrt{cf_0e}}.
    \end{bmatrix}$$
    By lower and upper bounds on $cf_0,e$ and upper bounds on $d$ we get $A_0$ satisfies the conclusion of the lemma. Let,
    $$\tilde{g} = g\circ A_0 = w_{\rho_m\circ A_0}, \; \tilde{h} = h \circ A_0,\; \tilde{v} = v\circ A_0.$$
    From equation \eqref{error_bound_lem1} we have that,
    $$|\tilde{h}-\tilde{g}| \leq C(t_0^{1+\frac{\delta}{2}}+\e_0^{\delta_m})$$
    for,
    $$\{(x,z)\;|\;|(x,z)|_m^2 \leq \frac{\beta}{C}t_0\} \subseteq A_0\begin{bmatrix}
        x\\z
    \end{bmatrix}\in S.$$
    By combining this with $|\tilde{h}-\tilde{v}|\leq C\sqrt{\e_0}$,
    $$|\tilde{v}-w_{\rho_m\circ A_0}| \leq C(e_0^{\frac{1}{2}}+\e_0^{\delta_m}+t_0^{1+\frac{\delta}{2}}) = C(\e_0,t_0).$$
    Choose the constants sufficiently small so that
    $$ t_0(1-\theta_0) \leq t_0-C(\e_0,t_0) \implies C(\e_0,t_0) \leq \theta_0 t_0.$$
    This implies,
    \begin{equation}
        \begin{split}
            S_{t_0(1-\theta_0)}&(w_{\rho_m\circ A_0},(0,0))\subseteq S_{t_0}(w_{\rho_m\circ A_0}+C(\e_0,t_0),(0,0))\subseteq\\
            &S_{t_0}(w_{\rho_m\circ A_0}-C(\e_0,t_0),(0,0)) \subseteq\{ |z|^{m+2}+|x|^2 \leq \frac{\beta}{C}t_0\}.
        \end{split}
    \end{equation}
     Thus,
    $$S_{t_0(1-\theta_0)}(w_{\rho_m\circ A_0},(0,0)) \subseteq \{\tilde{v} < t_0\}\subseteq S_{t_0(1+\theta_0)}(w_{\rho_m\circ A_0},(0,0))$$

\end{proof}

\begin{rmk}\label{rmk:exp_for_all_m}
    Note that from the proof of the preceding lemma it is immediate that the conclusion of the later Theorem \ref{ereg} holds when the right hand side is exactly $\rho_m$ for any $m>-1$ since the scale choice $t_0$ can be made arbitrarily small as there is no $\e_0$ term.
\end{rmk}

\begin{lemma}
    \label{OSI}
    Let $v$ in a section $S_1$ satisfy
    $$\det{D^2v} = cf\rho_m(z),$$
    $$S_{1-\theta}(w_m,(0,0)) \subseteq S_1 \subseteq S_{1+\theta}(w_m,(0,0)),$$
    for some constant $c>0$, and function $f(x,z)$. Suppose  $0<\theta\leq \theta_1$ and $0<\e\leq\min\{\e_1,\theta^8\}$, small enough such that if,
    $$|f-f_0|\leq \epsilon$$
    for some constant $f_0>0$. Then there exists $t_1>0$ and $\delta_0$ sufficiently small such that,
    $$AS_{t_1(1-\theta t_1^{\delta_0})}(w_{\rho_m\circ A},(0,0)) \subseteq S_{t_1} \subseteq AS_{t_1(1+\theta t_1^{\delta_0})}(w_{\rho_m\circ A},(0,0)),$$
    $$A =  \begin{bmatrix}
        a_{11} & a_{12} \\
        0 & a_{22} 
    \end{bmatrix}, \;\;\;\; |A-I| \leq C\theta,$$
    $$|cf_0-1| \leq C\theta,$$
    Here $\e_1,\theta_1,\delta_0$ depend upon $f_0$ and universal constants.
\end{lemma}
\begin{proof}
The proof is similar to \cite[Lemma 3.1]{SD} and our Lemma~\ref{BC}.
\end{proof}

\begin{theorem}\label{ereg}
    Let $v$ in a section $S_1$ satisfy
    $$\det{D^2v}(x,z) = f(x,z)\rho_m(z),$$
    $$S_{\gamma}(w_m,(0,0)) \subseteq S_1 \subseteq S_{\frac{1}{\gamma}}(w_m,(0,0)),$$
    for $f\in C^\alpha (S_1)$ and $\frac{1}{M}<f<M$. Now for some $\e$ small, if $[f]_{0,\alpha}=\sup\{|f-f(0,0)|/|x|^\alpha\} \leq \epsilon$, then for $m>-\beta_0$ (a positive constant $\beta_0$) there exists a matrix $A\in SL(2)$ such that $v$ has the following $2^{nd}$ order expansion,
    $$v(A(x,z)) =f(0,0)\int_0^z (z-t)(\rho\circ A)(t)dt + \frac{x^2}{2} + \bigO(|(x,z)|_m^{2+\delta}),$$
    $$A = \begin{bmatrix}
        a_{11} & a_{12} \\
        0 & a_{22} 
    \end{bmatrix}, \;\;\; c_1\leq|a_{ii}|\leq C_1,\;i=1,2 \;\text{ and }\; |a_{12}|\leq C_1.$$
    Here the error and constants $\e,c_1,C_1,\beta_0$ only depend upon, $M,\gamma$ and universal constants.
\end{theorem}
\begin{proof}
    We will run an Schauder type induction argument at the level of sections.\\
    
    \textit{Base Case:} From Lemma \ref{OSI} we have a constant $\theta_1$. Now let us choose $\theta_0 \leq \theta_1$ in Lemma \ref{BC} and thus choose $\e_0$ corresponding to this $\theta_0$. Now choosing $\e$ small such that $(diam(S)_1)^\alpha\e\leq \min\{e_0,e_1,\theta_0^8\}$ we get
    $$A_0S_{t_0(1-\theta_0)}(w_{\rho_m\circ A_0},(0,0)) \subseteq S_{t_0} \subseteq A_0S_{t_0(1+\theta_0)}(w_{\rho_m\circ A_0},(0,0)),$$
    $$A_0 =  \begin{bmatrix}
        \tilde{a}_{11} & \tilde{a}_{12} \\
        0 & \tilde{a}_{22} 
    \end{bmatrix}, \;\;\; c_0\leq|\tilde{a}_{ii}|\leq C_0,\;i=1,2 \;\text{ and }\; |\tilde{a}_{12}|\leq C_0.$$

    \textit{Induction Hypothesis:} for $k$, we have,
    $$A_0A_k S_{t_k(1-\theta_0 t_1^{(k-1)\delta})}(w_{\rho_m\circ A_0A_k},(0,0)) \subseteq S_{t_k} \subseteq A_0A_k S_{t_k(1+\theta_0 t_1^{(k-1)\delta})}(w_{\rho_m\circ A_0A_k},(0,0))$$
    
    \begin{equation}
        \label{MI}
    A_k = A_{k-1}B_{k-1},
    \end{equation}
    where $B_k$'s have the following structure,
    \begin{equation}
        \label{BI}
        B_k = \begin{bmatrix}
        b_{11,k} & b_{12,k} \\
        0 & b_{22,k} 
    \end{bmatrix},\; |b_{ii,k}-1|\leq C\theta_0 t_1^{(k-1)\delta},\;i=1,2\;\text{ and }\; |b_{12,k}|t_k^{-\frac{m}{2(2+m)}}\leq C\theta_0 t_1^{(k-1)\delta}.
    \end{equation}
    
    Here $t_k = t_0t_1^{k-1}$, $\delta \leq \min \{\delta_0,\frac{\alpha}{8}\min\{\frac{1}{2},\frac{1}{m+2}\}\}$, and the remaining constants are same as in Lemma~\ref{BC} and Lemma~\ref{OSI}. \\
    
    \textit{Induction Step:} Let
     $$\tilde{v}(x) = \frac{v(A_0A_k D_kx)}{t_k},$$
     where $D_k$ reflects the natural scaling
     $$D_k = \begin{bmatrix}
        t_k^{\frac{1}{2}} && 0\\
        0 && t_k^{\frac{1}{m+2}}
    \end{bmatrix}.$$
      Therefore $\{\tilde{v} < 1\} = D_k^{-1}A_k^{-1}A_0^{-1}S_{t_k}$ and
    $$S_{1-\theta_0 t_1^{(k-1)\delta}}(w_{t_k^{-\frac{m}{m+2}}\rho\circ t_k^{\frac{1}{m+2}}A_0A_k},(0,0))\subseteq \tilde{S}_1\subseteq S_{1+\theta_0 t_1^{(k-1)\delta}}(w_{t_k^{-\frac{m}{m+2}}\rho\circ t_k^{\frac{1}{m+2}}A_0A_k},(0,0)).$$
    Here $\tilde{v}$ solves the equation,
    \begin{equation}
        \begin{split}
                    \det(D^2\tilde{v})(x) = &t_k^{-\frac{m}{m+2}}\det(A_0A_k)^2f(A_0A_kD_kx)\rho_m((A_0A_kD_kx)_2)\\
                    & = c_kf_k(x)t_k^{-\frac{m}{m+2}}(\rho_m\circ t_k^{\frac{1}{m+2}}A_0A_k)(x) = c_kf_k(x)\rho_{m,k}(x).
            \end{split}
    \end{equation}    
    Further $f_k$'s satisfy, 
    \begin{equation}
        \begin{split}
        |f_k -f_0|\leq [f]_{0,\alpha}|A_0A_kD_kx|^{\alpha} = &[f]_{0,\alpha} |A_0A_kD_k|^{\alpha}|x|^\alpha\\
        & \leq C(diam(S_1))^\alpha\e (t_0^{\alpha}t_1^{(k-1)\alpha})^{\min\{\frac{1}{2},\frac{1}{m+2}\}}.        
        \end{split}
    \end{equation}
     By assumptions \eqref{MI} and \eqref{BI}, $|A_k|\leq C$. We once again choose $\epsilon$ small enough (Note that this choice is independent of the step we're at) to ensure that this difference is less than $\min\{\epsilon_1,(\theta_0t_1^{(k-1)\delta})^8\}$. Note that also, $\rho_{m,k}(x) = t_k^{-\frac{m}{m+2}}(\rho_m\circ t_k^{\frac{1}{m+2}}A_0A_k)(x)$ has similar constants as in inequality \eqref{dr}. Thus we can use Lemma \ref{OSI} to get
    $$\tilde{B}_kS_{t_1(1-\theta_0 t_1^{k\delta})} (w_{\rho_{m,k}\circ \tilde{B}_k},(0,0))\subseteq \tilde{S}_{t_1}\subseteq \tilde{B}_kS_{t_1(1+\theta_0 t_1^{k\delta})}(w_{\rho_{m,k}\circ \tilde{B}_k},(0,0)),$$
    Where $\tilde{B}_k$ is such that $|\tilde{B}_k-I|\leq C\theta_0t_1^{k\delta}$. Now to restate everything we need $B_k$ such that $B_kD_k = D_k\tilde{B}_k$. After a calculation we get our desired induction hypothesis in \eqref{BI} for $k+1$.
    Therefore,
    $$A_0A_{k+1} S_{t_{k+1}(1-\theta_0 t_1^{k\delta})} (w_{\rho_m\circ A_0A_{k+1}},(0,0))\subseteq S_{t_{k+1}} \subseteq A_0A_{k+1} S_{t_{k+1}(1+\theta_0 t_1^{k\delta})}(w_{\rho_m\circ A_0A_{k+1}},(0,0)).$$
   
    \textit{Finding a fixed coordinate change:} Note that $A^* = \lim_{k\to\infty} A_k$ exists and $A^* = A_kB_k^{*}$ where $B_k^* = \prod_{i\geq k}B_k$. This convergence is possible as long as $m>-2\delta/(1+\delta)$. Thus $A_k = A^*(B^*_k)^{-1} $. Using our natural scaling $D_k$,
    $$S_{t_k(1-\theta_0 t_1^{(k-1)\delta})}(w_{\rho_m\circ A_0A_k},(0,0)) = D_kS_{1-\theta_0 t_1^{(k-1)\delta}}(w_{t_k^{-\frac{m}{m+2}}\rho_m\circ t_k^{\frac{1}{m+2}}A_0A_k},(0,0)).$$
    We now find a matrix $C_k$ such that $(B^*_k)^{-1}D_k = D_kC_k$. After a direct calculation, $|C_k-I|\leq C\theta_0t_1^{(k-1)\delta}$. 
    Therefore we have the following chain of containments,
    \begin{equation}
        \begin{split}
           & D_kS_{1-\theta_0 t_1^{(k-1)\delta}}(w_{t_k^{-\frac{m}{m+2}}\rho_m\circ t_k^{\frac{1}{m+2}}A_0A_k},(0,0))\\
            \subseteq \; &B_k^*D_kC_kS_{1-\theta_0 t_1^{(k-1)\delta}}(w_{t_k^{-\frac{m}{m+2}}\rho_m\circ t_k^{\frac{1}{m+2}}A_0A_k},(0,0))\\
            \subseteq \;& B_k^*D_k S_{1- C\theta_0 t_1^{(k-1)\delta}}(w_{t_k^{-\frac{m}{m+2}}\rho_m\circ t_k^{\frac{1}{m+2}}A_0A_kC_k^{-1}},(0,0))\\
            \subseteq \;& B_k^*S_{t_k(1-C\theta_0 t_1^{(k-1)\delta})}(w_{\rho_m\circ A_0A^*(B_k^*)^{-1}C_k^{-1}},(0,0)).\\
        \end{split}
    \end{equation}
    We have used the fact that $w_\rho\circ C_k^{-1} = w_{\rho\circ C_k^{-1}}+\bigO(t_k^{\delta})$.
    Now since all our matrices are upper triangular and the only coordinate we need for the $\rho$ term is $z$, we have,
    $$\rho_m\circ A_0A^*(B_k^*)^{-1}C_k^{-1} = \rho_m\circ A_0A^*D_kC_kD_k^{-1}C_k^{-1} = \rho_m\circ A_0A^*.$$
    Therefore,
    $$A_0A^* S_{t_{k+1}(1-C\theta_0 t_1^{k\delta})} (w_{\rho_m\circ A_0A^*},(0,0))\subseteq S_{t_{k+1}} \subseteq A_0A^* S_{t_{k+1}(1+C\theta_0 t_1^{k\delta})}(w_{\rho_m\circ A_0A^*},(0,0)),$$
    and we can conclude our second order expansion. Note that $A\in SL(2)$ follows from using the equation.

\end{proof}

\begin{rmk}
    In the above proof we only needed a uniform $C^{\alpha}$ estimate in the $x$ variable since we can rewrite the right hand side as,
    $$f(x,z)\rho_m(z) = \frac{f(x,z)}{f(0,z)}f(0,z)\rho_m(z)$$
     Then by defining our new $\tilde{f}(x,z) = \frac{f(x,z)}{f(0,z)}$, $[\tilde{f}]_\alpha$ is, 
     $$\frac{|\tilde{f}-\tilde{f}_0|}{|(x,z)-(0,0)|^\alpha} \leq M \frac{|f(x,z)- f(0,z)|}{|x|^{\alpha}} \leq M[f]_{x,\alpha}.$$
     Therefore the theorem also holds under incomplete H$\ddot{o}$lder conditions as in \cite{Li}.
\end{rmk}

\begin{corollary}\label{expansion_corollary}
    Let $v$ in a section $B_c \subseteq S_1\subseteq B_C$ satisfy
    $$\det{D^2v}(x,z) = f(x,z)\rho_m(z)$$
    for $f\in C^\alpha (S_1)$, $\frac{1}{M}<f<M$ and $m>-\beta_0$. Then there exists a matrix $A$ such that $v$ has the following $2^{nd}$ order expansion,
    $$v(A(x,z)) = f_0\int_0^z (z-t)(\rho\circ A)(t)dt + \frac{x^2}{2} + \bigO(|(x,z)|_m^{2+\delta}),$$
    $$A = \begin{bmatrix}
        a_{11} & a_{12} \\
        0 & a_{22} 
    \end{bmatrix}, \;\;\; c_1\leq|a_{ii}|\leq C_1,\;i=1,2 \;\text{ and }\; |a_{12}|\leq C_1.$$
    Here the error and constants $\beta_0,c_1,C_1$ only depend upon, $[f]_{\alpha},M$ and universal constants.
\end{corollary}

\begin{proof}
    By the balancing Lemma \ref{balacing}, $A_tB_c \subseteq S_t\subseteq A_tB_C$, with
    $$A_t = c_t\begin{bmatrix}
        \lambda_t && b_t \\
        0 && \frac{1}{\lambda_t}
    \end{bmatrix}.$$
     Also note that by Lemma \ref{SCR},
    $$L_{x} + \omega(|x-y|) \leq u(y)$$
    for a universal modulus of strict convexity $\omega$ independent of a point. This shows that $(B_{\omega^{-1}(h)})^c \subseteq (S_h(x))^c$. Thus, $A_hB_c \subseteq S_h \subseteq B_{\omega^{-1}(h)}$. Thus $|A_h| \leq \omega^{-1}(h)$. Therefore by zooming in till $h=t_0$ and considering,
    $$v_h = \frac{v(A_hx)}{h}.$$
    $v_h$ satisfies the equation
    $$\det(D^2v_h) = \frac{c_h^4}{h^2}f(A_hx)\rho_m((A_hx)_2) = C_h f_h \rho_m(\frac{c_h}{\lambda_h}z)$$
    with $B_c \subseteq S_1 \subseteq B_C$ and $f_h$, $[f_h]_\alpha \leq |A_h|^\alpha[f]_\alpha$. By choosing $t_0$ small but fixed depending upon $[f]_{\alpha}$, $v_h$ satisfies the conditions of Theorem \ref{ereg} and the expansion of $v$ follows.
\end{proof}

The following theorem follows immediately from Corollary \ref{expansion_corollary}.
\begin{theorem}\label{C^2_estimate}
    Let $v$ in a section $B_c\subseteq S_1\subseteq B_C$ satisfy
    $$\det{D^2v}(x,z) = f(x,z)\rho_0(z)$$
    for $f\in C^\alpha (S_1)$, $\frac{1}{M}<f<M$. Then,
    $$\|D^2v\|_{L^{\infty}(S_{1/2})} \leq C.$$
    Here the constant $C$ depends on $[f]_\alpha,M$ and universal constants.
\end{theorem} 

\begin{rmk}
    For $m>0$, a similar $C^{1,1}$ estimate follows immediately by using standard arguments. See for e.g. Savin-Le in \cite{SL}. For $m<0$ we can bound $\nabla v_x$.
\end{rmk}

\section{A Liouville Theorem for $\det(D^2u) = \rho_m(z)$}\label{sec:Liouville}

Since we have obtained a $2^{nd}$ order expansion for solutions to the Monge-Ampere equation \eqref{doubling_equation}, one might expect global solutions to also have rigidity. We show that it is indeed true. We still assume $m>-1$, $S_t = S_t(u,(0,0))$, $u(0)=0$ and $\nabla u(0,0) = 0$. We have the following important balancing of section.
\begin{lemma}\label{bcing}
    Let $u$ be a strictly convex solution to,
    \begin{equation}
    \begin{split}
        \det(D^2u) =& \rho_m(z)\;\;  \text{ on }\;\; S_1,\\
         \end{split}
    \end{equation}
    such that $u_{xx}(0) = 1$ and $u_{xz}(0) = 0$. Then for some universal $\beta>0$, the section $S_1$ is balanced as follows,
    $$B_\beta \subseteq S_1 \subseteq B_{\frac{1}{\beta}}.$$
\end{lemma}
\begin{proof}
    Consider the normalization $\tilde{u}$,
    $$\tilde{u}(x) = \frac{1}{\lambda}(u(T^{-1}x)-1),$$
    which satisfies the equation
    $$\det(D^2\tilde{u}) = \frac{\rho_m((T^{-1}x)_2)}{(\lambda\det(T))^2}= \tilde{\rho}_m(z).$$
    Here $T$ is such that $B_1\subseteq T(S_1)\subseteq B_n$, and $\lambda$ is chosen such that the total Monge-Ampere mass, $\mathcal{M}(\tilde{u})(T(S_1))=1$. We assume,
    $$T = \begin{bmatrix}
        t_{11} && t_{12} \\
        0 && t_{22}
    \end{bmatrix}.$$
    From \cite{CG1} the normalized solutions satisfy
     $$\tilde{C}_1\leq |\min_{T(S_1)}\tilde{u}|\leq \tilde{C}_2,$$
      where these constants only depend upon universal constants and the doubling constant of $\rho_m$. Since $\min_{x\in S_1}\{u(x)\} = 0$, $\lambda$ has to be bounded between the above constants. Moreover,
     \begin{equation}\label{bl}
        \iint_{T(S_1)} \tilde{\rho}_m(z)dxdz = 1.
    \end{equation}
     Thus $\tilde{c}_1|z|^m\leq \tilde{\rho}(z)\leq \tilde{c}_2|z|^m$ where $\tilde{c}_1,\tilde{c}_2$ also depend on universal constants. Therefore Theorem \ref{ereg} shows that,
     $$\tilde{u}(Ax) =\min \{\tilde{u}\}+ w_{\tilde{\rho}_m}+\bigO(|(x,z)|_m^{2+\delta}).$$
     Thus $\nabla \tilde{u}_x = A^{-tr}diag(1,0)A^{-1} \textbf{e}_1$ where $\textbf{e}_1 = (1,0)^{T}$. However,
     $$\nabla\tilde{u}_x(0) = \frac{1}{\lambda}(T^{-1})^{tr}diag(1,0)T^{-1}\textbf{e}_1$$
     by definition. Equating both, we get appropriate bounds on $t_{11},t_{12}$. Then using equation \eqref{bl} we get double sided estimates for $t_{22}$ as well. This gives us balancing.
\end{proof}

We now show that global solutions to the equation \eqref{doubling_equation} have to be strictly convex,
\begin{lemma}\label{no_line}
        Let $u$ be a global solution of $\det(D^2u) = \rho_m(z)$ on $\mathbb{R}^2$, then $u$ is strictly convex.
    \end{lemma}
    \begin{proof}
        We prove by contradiction. For $z\neq 0$, since $\det(D^2u)>0$, $u$ has to be strictly convex. Thus singularities can only occur along $\{z=0\}$. However by propagation of singularities, if $u$ is not strictly convex, the singularity should propagate. Thus $u$ must be singular along entire line $\{z=0\}$. By subtracting an appropriate linear function, $u=0$ on $\{z=0\}$. Now for any $p=(x,z)$ we show that $\partial(u)(p) \subseteq \{0\}\times \mathbb{R}$. Assume not, then let $(a,b)\in \partial(u)(p)$ with $a\neq 0$. Therefore on $z=0$, the line $L(q,0) = u(p)-a(q-x)-b(0-z)$ must be bounded above by $0$. However if $a\neq 0$ we can choose $q$ large such that we contradict this. Thus for any $p$, $\partial(u)(p)\subseteq \{0\}\times \mathbb{R}$. Thus $\mathcal{M}u(\mathbb{R}^2)= 0$, a contradiction. 
    \end{proof}

Now that strict convexity is guaranteed, we prove rigidity by applying Lemma \ref{bcing} to appropriately scaled sections of $u$. This proves Theorem \ref{Liouville1}.
\begin{theorem}\label{liouville_1}
    Let $u$ be a solution to,
    \begin{equation}\label{liouville_pde}
    \begin{split}
        \det(D^2u) = \rho_m(z)\;\;  \text{ on }\;\; \mathbb{R}^2\\
         \end{split}
    \end{equation}
    Then up to a fixed coordinate change denoted by matrix $A$ and some constant $a>0$, 
    $$u(Ax) =u = \frac{a}{2}x^2 + \frac{1}{a}\int_0^{z}\rho(t)(z-t)\,dt.$$
\end{theorem}
\begin{proof}
    We assume $u(0)=0$, $\nabla u(0)=0$ and $u\geq 0$ by subtracting the supporting hyperplane at $0$. Moreover by an upper triangular coordinate change, we also ensure that $\nabla u_x(0) = (1,0)^T$ since $u$ has an second order expansion as in Corollary \ref{expansion_corollary} (for $-1<m<0$ see Remark \ref{rmk:exp_for_all_m}). For a section $S_h$ normalize it by considering $u_h(x) = \frac{1}{h}u(D_hx)$, where,
    $$D_h = \begin{bmatrix}
        h^{\frac{1}{2}} && 0\\
        0 && h^{\frac{1}{m+2}}
    \end{bmatrix}.$$
    Thus $u_h$ satisfies,
    \begin{equation}
    \begin{split}
        \det(D^2u_h) = \rho_{h,m}&(z) = h^{-\frac{m}{m+2}}\rho_m(h^{\frac{1}{m+2}}z)\;\;  \text{ on }\;\; S_1\\
         \end{split}
    \end{equation}
    By Lemma \ref{bcing}, $S_1$ is balanced around $0$ and Theorem \ref{ereg} gives,
    $$u_h\circ A_h(x,z) = \frac{x^2}{2}+\int_0^z(z-t)\rho_{h,m}((A^h)_{22}t)dt + \bigO(|x,z|_m^{2+\delta}).$$
    Rewriting,
    $$u(x,z) = \frac{(a_{11}^hx+h^{\frac{m}{m+2}}a_{12}^hz)^2}{2}+\int_0^{a_{22}^hz}(a_{22}^hz-t) \rho_m (\frac{t}{a^h_{22}})dt+\frac{\bigO(|(x,z)|_m^{2+\delta})}{h^{\delta/2}}.$$
    Here $[A_h^{-1}] = a^h_{ij}$ for $i,j=1,2$. Since $u_{xx}=1 $ and $u_{xz}=0$ at origin we get that the mixed term $a_{12}^h$ is identically $0$. Note that the family $A_h$ is compact. Thus taking a subsequence $h_n\to\infty$ such that $A_{h_n}\to A_0$. Therefore for a fixed $(x,z)$,
    $$u(x,z) = \frac{(a_{11}x)^2}{2}+a_{22}^2\int_0^z (z-t)\rho_m (t)dt.$$
    The relation between $a_{11},a_{22}$ follows from the equation \eqref{liouville_pde}.
\end{proof}

\section{Regularity at the minimum point}\label{sec:min_point_analysis}

In this section we are concerned with the following theorem,

\begin{theorem}\label{x=0_theorem}
    Let $u$ be a convex solution to,
    \begin{equation}\label{x=0_pde}
        \det(D^2u) = f_p(u)
    \end{equation}
    in a section $B_c \subseteq S_1(u,(0,0))\subseteq B_C$ such that $u(0) =0$,$ \nabla u(0,0) = (0,0)$ is the isolated minima. Here 
    $$c_0t^p \leq f_p(t) \leq C_0t^p$$
    for $-2< p<2$ and $c_0,C_0$ constants. Then for every $t \leq 1$,
    $$B_{c_1\sqrt{t^{1-p/2}}} \subseteq S_t(u,(0,0)) \subseteq B_{C_1\sqrt{t^{1-p/2}}}$$
    for constants $c_1,C_1$ depending upon $c,C,c_0,C_0,p$. Moreover, up to a universally bounded continuous affine change $L(t)\in ASL(2)$, $r(t)^{-1}L(t)S_t \to B_1(0)$. Here $r(t) = \sqrt{\frac{A(t)}{\pi}}$ where $A(t)$ is the area of the section $S_t$.
\end{theorem}

In the theorem above $ASL(2) = SL(2)\rtimes \mathbb{R}^2$, i.e., the group generated by translations and unimodular transforms. There has been a great deal of interest in trying to study the level sets of solutions to \eqref{x=0_pde}. Assuming the above theorem is true we give a complete proof of Theorem \ref{Main}.

\begin{proof}[Proof of Theorem \ref{Main}]
    To prove \Romannum{1}, we use Theorem \ref{C^2_estimate} and combine them with bounded rescaling using Theorem \ref{x=0_theorem}. For $T_t = r(t)^{-1}L(t)$, let $u_t = \frac{1}{t}u(T_t^{-1}x)$. Therefore $u_t$ satisfies,
    $$\det(D^2u_t) = \tilde{f}(u_t) \;\text{ on } B_{1/2}\subseteq T_tS_t\subseteq B_{2}$$
    for $t\leq t_0$ with $t_0$ universal. By Theorem \ref{C^2_estimate}, 
    $$\|D^2u_t\|_{L^{\infty}(S_{1/2})} \leq C.$$
    Since $L(t)$ are universally bounded, we immediately get
    $$\|D^2u\|_{L^{\infty}(S_{t}\setminus S_{t/2})} \leq C$$
    for all $t\leq t_0$. Moreover, at $t_0$, we have directly have $\|D^2u\|_{L^{\infty}(S_{1/2}\setminus S_{t_0/2})} \leq C$ from Theorem \ref{C^2_estimate}. Therefore we get \Romannum{1}.\\
    
    The smoothness claimed in \Romannum{2} follows from higher order estimates done in \cite{Li} and the expansion follows from Corollary \ref{expansion_corollary}.\\
    
    \Romannum{3} follows immediately from Theorem \ref{x=0_theorem}.
\end{proof}

\begin{rmk}\label{rmk:gen}
    Note that since we only use Theorems \ref{C^2_estimate} and \ref{x=0_theorem}, the statement of Theorem \ref{Main} also holds if $f$ is allowed to touch zero at finitely many points in a polynomial fashion (equation \eqref{dr}) as discussed in Section \ref{subsec:PLT}.
\end{rmk}

\subsection{$2^{nd}$ order flow and Monotone quantity}\label{sec:flow}
We cast the problem into a second order flow for the support function of the level sets. Let $h(t,\theta)$ be the support function of $K_t = \{u<t\}$. Define $v = h_t$. The following relations hold
\begin{equation}
    Du = \frac{1}{v}n(\theta),\;\;\rho = h+h_{\theta\theta}.
\end{equation}
Here $\rho$ is the radius of curvature for the level set $K_t$ and $n(\theta) = (\cos(\theta),\sin(\theta))$. Also note $v = 1/|Du|>0$. The area $A(t)$ of a level set $K_t$ can be written as
$$A(t) = \frac{1}{2}\ic h^2-h_\theta^2d\theta.$$
Therefore differentiating in $t$ and integrating by parts,
$$A'(t) = \int_{\mathbb{S}^1}\rho vd\theta$$
and 
\begin{equation}\label{double_der_A}
    A''(t) = \int_{\mathbb{S}^1} [v^2 - v_{\theta}^2]\,d\theta + \int_{\mathbb{S}^1} \rho v_t d\theta.
\end{equation}
In the sequel we will be interested in capturing the non-homothetic changes in $K_t$ as $t$ varies. To do so we use the Alexandrov-Fenchel inequality (see \cite{Sch}), which implies that
$$d(s) := V(K_{t+s},\,K_t)^2 - A(t+s)A(t) \geq 0, \quad d(0) = 0,$$
where the mixed volume in terms of support functions is given by
$$V(K_{t+s},\,K_t) = \frac{1}{2}\int_{\mathbb{S}^1} [h(t+s,\,\theta)h(t,\,\theta) - h_{\theta}(t+s,\,\theta)h_{\theta}(t,\,\theta)]\,d\theta,$$
and that equality holds if and only if $K_{t+s}$ is homothetic to $K_t$. Differentiating twice and using the above identities for $A$ and its derivatives we get
\begin{equation}\label{AF}
0 \leq d''(0) = \frac{1}{2}A'^2(t) - A\int_{\mathbb{S}^1} [v^2-v_{\theta}^2]\,d\theta.
\end{equation}

We now turn our attention to the PDE \eqref{x=0_pde}. Writing $f_p=f$ and using notations introduced above, we get
\begin{equation}\label{2nd_order_evol}
    v_t = -f(t)\rho v^4 \iff (v^{-3})_t = 3f(t)\rho.
\end{equation}
Using the above PDE we can write,
\begin{equation}\label{ga}
    G' = 2f(t)A'. 
\end{equation}
where
\begin{equation}\label{definition_G}
    G = \ic v^{-2} = \ic |Du|^2.
\end{equation}
The equation \eqref{ga} is the Euler-Lagrange equation of the functional
\begin{equation}\label{EL}
 \int_{t_0}^{t_1} \left[\ic v^{-2}d\theta-6f(t)A(t)\right] dt = \int_{t_0}^{t_1}G(t)-6f(t)A(t)dt.\end{equation}

By the Cauchy-Schwarz inequality,
\begin{equation}\label{MPQuant}
(A')^2 \leq \left(\int_{\mathbb{S}^1} \rho^2v^4\,d\theta\right)G,
\end{equation}
with equality if and only if $\rho v^3$ is constant. 

\begin{rmk}
It was noted in \cite{CMS}, \cite{FS} that the quantity $\rho v^3$ takes its maximum on the boundary when $f$ is increasing, which allowed the authors in those references to control the curvature of level sets. 
\end{rmk}

 Combining the observations \eqref{AF} and \eqref{MPQuant} we have the following monotone quantity,
\begin{lemma}\label{monotone}
    For a solution to the PDE \eqref{x=0_pde}, the quantity
    $$Q = G\frac{A'^2}{A} $$
    is monotone. Moreover
    \begin{equation}\label{dissipation}
    \frac{d}{dt}\log(Q) = \frac{-2}{A'}(\delta_1+f\delta_2) \leq 0
\end{equation}
for the defects $\delta_1,\delta_2$ defined as,
\begin{equation}\label{defect_definition}
    -\delta_1 = \int_{\mathbb{S}^1}v^2-v_\theta^2d\theta-\frac{(A')^2}{2A}\leq 0,\;\;\;-\delta_2 = \frac{(A')^2}{G} - \ic \rho^2v^4\leq 0
\end{equation}
\end{lemma}
\begin{proof}
    Differentiate $\log(Q)$ and use \eqref{double_der_A}, \eqref{AF} and \eqref{MPQuant}.
\end{proof}

\begin{rmk}
    In the above dissipation relation we note that from \eqref{AF} that $\delta_1$ governs homothety. Therefore $\delta_1=0$ implies sections are evolving by homothety and $v = a(t)h$ for some function $a(t)$ up to translations. As mentioned before $\delta_2=0$ implies $\rho v^3 = c(t)$ for some constant depending only on the height $t$. Therefore since $\rho$ is translation invariant we get an equation of the form $\rho h^3 = \tilde{c}(t)$. We note that the only solutions to the equation are support functions of ellipses (see Lemma \ref{normalization}).
\end{rmk}

We have the following important observation regarding the defects and relevant quantities. By affine invariant, we mean invariant under the action of the group $ASL(2)$.
\begin{lemma}\label{aff_inv}
    The defects $\delta_1,\delta_2$ and the quantities $G,\ic[ v^2-v_{\theta}^2]d\theta, \ic \rho^2v^4d\theta$ are affine invariant.
\end{lemma}
\begin{proof}
    Since $\delta_1$ comes from an Alexandrov-Fenchel defect, it has to be affine invariant. Since 
    $$\ic v^2-v_{\theta}^2d\theta = -\delta_1+\frac{(A')^2}{2A},$$
    it has to be affine invariant. By \eqref{ga}, $G$ has to be affine invariant. By \eqref{double_der_A}, $\ic \rho^2v^4d\theta$ has to be affine invariant. Thus $\delta_2$ has to be affine invariant by definition.
    
\end{proof}

\subsection{Rescaled Flow}\label{rescale}
We now perform an area normalized rescaled flow argument by forcing the area of every $K_t$ to be $\pi$. Set
$$
r(t):=\sqrt{\frac{A(t)}\pi},\qquad
\tau:=\log r,\qquad
h(t,\theta)=r(t)H(\tau,\theta).
$$
Then the body with support function $H$ has area $\pi$. Define
$$
R:=H+H_{\theta\theta},\qquad
Z:=H_\tau,\qquad
W:=H+Z.
$$
Here $R$ is the radius of curvature for the convex body associated with $H$. We also call $Z$ the physical velocity. Since
$$
\tau_t=\frac{A'}{2A},\qquad r_t=r\tau_t,
$$
we get
\begin{equation}
v=h_t=r_tH+rH_\tau\tau_t
=\frac{A'}{2\sqrt{\pi A}}W.
\end{equation}
Write
$$
\beta:=\frac{A'}{2\sqrt{\pi A}},\qquad v=\beta W.
$$
The area identity becomes
\begin{equation}\label{fix_area}
\ic(H^2-H_\theta^2)\,d\theta=2\pi.
\end{equation}
Differentiating \eqref{fix_area} in $\tau$ gives
\begin{equation}\label{orthogonality}
\ic(HZ-H_\theta Z_\theta)\,d\theta=0.
\end{equation}
Also, $\rho=rR$, so
$$
A'=\int\rho v
=r\beta\int RW
=\frac{A'}{2\pi}\int RW,
$$
and hence
\begin{equation}
\ic RW\,d\theta=2\pi.
\end{equation}
We now rewrite everything in terms of the newly defined quantities. If we define,
$$I = \ic W^{-2}d\theta,$$
 Since $v=\beta W$,
$$
G=\beta^{-2}I=\frac{4\pi A}{A'^2}I,
$$
and
\begin{equation}\label{MQ}
Q=\frac{GA'^2}{A}=4\pi I.
\end{equation}
We now first convert $\delta_1$. Note that $2\pi\beta^2=A'^2/(2A)$, therefore we have,
$$-\delta_1 = \beta^{2}(\beta^{-2}\int_{\mathbb{S}^1}v^2-v_\theta^2d\theta-2\pi) =\beta^{2}(\ic W^2-W_\theta^2d\theta-2\pi).$$
Next,
\begin{equation}
    \begin{split}
        \int(W^2-W_\theta^2)
&=\int(H^2-H_\theta^2)
+2\int(HZ-H_\theta Z_\theta)
+\int(Z^2-Z_\theta^2)\\
& = 2\pi + \ic(Z^2-Z_\theta^2)d\theta
    \end{split}
\end{equation}
Therefore, if we define
 $$E_1 = \ic(Z_\theta^2-Z^2)d\theta,$$
 we get $\delta_1 = \beta^2E_1$. Also
$$
\ic\rho^2v^4d\theta=r^2\beta^4\int R^2W^4
=\frac{A'^4}{16\pi^3A}\int R^2W^4,
$$
while
$$
\frac{A'^2}{G}
=\frac{A'^4}{4\pi AI}
=\frac{A'^4}{16\pi^3A}\frac{4\pi^2}{I}.
$$
Thus define
$$E_2=\int R^2W^4\,d\theta-\frac{4\pi^2}{I}.$$
We get $\delta_2=\frac{A'^4}{16\pi^3A}E_2 = \frac{A'^2\beta^2}{4\pi^2}E_2$. Note that we can rewrite $E_2$ as,
$$
\begin{aligned}
\int\left(RW^2-\frac{2\pi}{I}W^{-1}\right)^2d\theta
&=\int R^2W^4
-\frac{4\pi}{I}\int RW
+\frac{4\pi^2}{I^2}\int W^{-2}\\
&=\int R^2W^4-\frac{8\pi^2}{I}+\frac{4\pi^2}{I}\\
&=E_2.
\end{aligned}
$$
Since $d\tau/dt=A'/(2A)$, \eqref{dissipation} gives,
\begin{equation}\label{dissipation_rel}
-\frac d{d\tau}\log I
=\frac1\pi\bigl(E_1+a(\tau)E_2\bigr)\geq 0,
\qquad
a(\tau):=\frac{f(t)A'(t)^2}{4\pi^2}.
\end{equation}
As long as $a(\tau) \sim 1$ both the terms contribute to the dissipation and if $Q \sim 1$, we have,
\begin{equation}\label{finite_dissipation}
    \int_{-\infty}^T (E_1+E_2) \leq M <\infty
\end{equation}
for some $M$ depending on universal constants. We show that under the conditions of Theorem \ref{x=0_theorem} the necessary quantities are indeed universally bounded from above and below.
\begin{lemma}\label{rescaling}
Let $f_p$ satisfy the conditions of Theorem \ref{x=0_theorem}. Let $u$ be a solution to \eqref{x=0_pde} with an isolated minimum,
$$\min \{u\} = u(0) = 0.$$
Let $A,G,Q$ be quantities introduced in Section \ref{sec:flow}. Then we have the following bounds,
$$A(t) \sim t^{1-\frac{p}{2}},\qquad A'(t)\sim t^{-p/2},\qquad
G(t)\sim t^{1+p/2},\qquad
Q(t)\sim1.$$
The constants depend only upon universal constants and $c_0,C_0$. 
    
\end{lemma}
\begin{proof}
    For sections $S_t = S_t(u,(0,0))$, normalize the sections,
    $$B_{c_J}\subseteq T_t S_t \subseteq B_{C_J}.$$
    Then $u_t = t^{-1}u(T_t^{-1}(x-c_t))$ satisfies,
    $$\det(D^2u_t) = t^{-2}\det(T_t)^{-2}f_p(tu_t) = f_t(u_t).$$
    If $c_0\kappa_t t^p \leq f_t(u_t) \leq C_0 \kappa_t t^p$ then we show $\kappa_t \sim 1$. Here $\kappa_t = \det(T_t)^{-2}t^{p-2}$. 
    For the upper bound consider the barrier function $\Phi$,
    $$
    \Phi(y)
    =
    \left(\frac{(|y|-r_0)_+}{c_J-r_0}\right)^\alpha.
$$
Here $r_0 = c_J/2$ and $\alpha > \max\left\{2,\frac3{2-p}\right\}$. Note that $\Phi$ satisfies,
$$
\Phi\ge1\quad\hbox{on }\partial T_tS_t,
$$
$$
\Phi\equiv0\quad\hbox{on }B_{c_J/2},
$$
and
\begin{equation}
\det D^2\Phi\le\min\{ C_B\Phi^p,C_B'\}
\quad\hbox{in }T_tS_t.
\end{equation}
This now forces $c_0\kappa_t \leq \max\{C_B,C_B'\}$. If not then by comparison principle, $u_t \leq \Phi$ in $T_tS_t$. Since $u_t =1 \leq \Phi$ on $\partial T_tS_t$, then if  $u_t>\Phi$ in the interior, then for some $y\in T_tS_t$,
$$u_t(y) > \Phi(y),\qquad D^2u_t \leq D^2\Phi. $$
 However
 $$\det(D^2u_t) \geq c_0\kappa_t u_t^p > \min\{C_B \phi^p,C_B'\} \geq \det(D^2\Phi),$$
 a contradiction. This is not possible since $u_t = 0$ only at one point. Therefore $\kappa_t \leq C$ and 
 \begin{equation}\label{A_upperbound}
     A(t) \leq C  t^{1-p/2}.
 \end{equation}
 Now for the lower bound, the Alexandrov maximum principle gives,
    $$|\min (u_t-1)| = 1 \leq C\left(\int_{T_tS_t}C_0\kappa_t u_t^p \right)^{\frac{1}{2}}.$$
    Note that if $p\geq 0$ then this immediately gives a lower bound for $\kappa_t$. If $p<0$, by Cavalieri's principle we have
    $$\int_{T_tS_t} u_t^p = \int_0^1 s^pdA_{u_t}(s) = A_{u_t}(1) - p\int_0^1s^{p-1}A_{u_t}(s)ds.$$
    Note that,
    $$A_{u_t}(s) = |\{u_t = s\}| = |T_t\{u=st\}|\leq \kappa_t^{-\frac{1}{2}}t^{p/2-1}A(st)$$
    By \eqref{A_upperbound} we can bound $s^{p-1}A_{u_t}(s)$ in terms of $\kappa_t^{-1/2} s^{p/2}$. Therefore,
    $$\int_{T_tS_t} u_t^p \leq C(1+\kappa_t^{-1/2})$$
    as long as $s^{p/2}$ is integrable, i.e., $p>-2$. Therefore we get a lower bound for $\kappa_t$ even in this case. Thus,
 $$A(t) \sim t^{1-p/2}.$$
 By Brunn-Minkowski theory $\sqrt{A(t)}$ has to be concave and we get
 $$A'(t) \sim t^{-p/2}.$$
 The bounds for $G,Q$ follow from \eqref{ga}.

\end{proof}

\begin{rmk}
    The Alexandrov maximum principle step above forces that no isolated minima can occur when $p \geq 2$ since 
    $$\frac{1}{C} \leq \kappa_t \leq \frac{1}{c_J^2}A(t)^2t^{p-2}. $$
    For $p>2$, this forces $A(t)\to \infty$ as $t\to 0$ which isn't possible.
\end{rmk}

\begin{rmk}\label{rmk:large_t}
    Note that when $K = \{u=0\}$ is a set of finite measure then the conclusion of the above Lemma still holds for $t$ large enough. This is because from the Alexandrov maximum principle step we get a lower bound for the area $A(t)$. Therefore for a large $t$, the measure of the set $\{u=0\}$ can be made arbitrarily small after normalizing by choosing $t$ sufficiently large. Then the barrier argument forces $T_tK \supset B_{r_0}$ for a fixed constant $r_0$, which contradicts the small measure of $T_tK$.
\end{rmk}

We denote,
$$\mathcal{E}(\tau) = E_1(\tau)+E_2(\tau). $$
We note again that the dissipations as well as the quantities $I$ are affine invariant. Moreover, since area $=\ic h^2-h_\theta^2d\theta$ is affine invariant, therefore the bilinear form $\ic h_1h_2-(h_1)_\theta(h_2)_\theta$ is also affine invariant.

\subsection{Compactness for Affine Normalized Sections} In this section we show a compactness result for small dissipations. We first note the following affine normalized dissipation inequality,

\begin{lemma}\label{Z_error}
Let $(H(\tau,\cdot),Z(\tau,\cdot))$ be the support function and the physical velocity of fixed area convex bodies $B_{c_1} \subseteq K_{\tau}^* \subseteq B_{C_1}$. Here $c_1,C_1$ are some fixed constants. Then,
$$
\|Z-P_1Z\|_{H^1}^2\leq CE_1,
$$
with $C$ depending only on $c_1,C_1$.
\end{lemma}

\begin{proof}
If $H_k = P_k(H)$ denote the Fourier modes, then for
$$
S_H:=\sum_{k\ge2}(k^2-1)|H_k|^2,
$$
the fixed area condition \eqref{fix_area} gives
$$
H_0^2-S_H=1.
$$
The condition $B_{c_1} \subseteq K_{\tau}^* \subseteq B_{C_1}$ gives  $c_1 \leq H_0 = \frac{1}{2\pi}\ic H d\theta \leq C_1$. So
$$
\frac{S_H}{H_0^2}\le1-c_J
$$
for some $c_J>0$ depending upon $c_1,C_1$. The orthogonality \eqref{orthogonality} expresses $Z_0$ as
$$
H_0Z_0=\sum_{k\ge2}(k^2-1)H_k\cdot Z_k.
$$
Cauchy--Schwarz therefore gives
$$
|Z_0|^2\le(1-c_J)\sum_{k\ge2}(k^2-1)|Z_k|^2.
$$
Since
$$
E_1=\sum_{k\ge2}(k^2-1)|Z_k|^2-|Z_0|^2,
$$
we obtain a uniform positive fraction of the $H^1$ norm of all modes except $k=1$.
\end{proof}

\begin{lemma}\label{normalization}
    Let $K_t$ be convex bodies evolving by \eqref{2nd_order_evol}. Let $(H(\tau,\cdot),Z(\tau,\cdot))$ be the support function and the physical velocity of fixed area convex bodies that have been affine normalized such that $B_{c_1} \subseteq K_{\tau}^* \subseteq B_{C_1}$. Here $K_{\tau}^*$ corresponds to the section $K_t$.  Given $\epsilon>0$, there exists $\delta>0$ such that if $\mathcal{E}(\tau) \leq \delta$, then there exists $L(\tau)\in ASL(2)$ with a symmetric $SL(2)$ part such that for $(\tilde{H},\tilde{Z})$ associated to $L(\tau) K_\tau^*$  we have,
    \begin{equation}\label{close_circle}
        \|\tilde{H}-1\|_{H^2}+\|\tilde{Z}\|_{H^1} \leq \epsilon
    \end{equation}
\end{lemma}
\begin{proof}
    Assume not. Then for a given $\epsilon>0$, for every $j\in \mathbb{N}$, there exist times $\tau_j$ and associated functions $(H_j,Z_j)$ such that $\mathcal{E}(\tau_j) \leq 1/j$ but 
    $$\|\tilde{H}_j-1\|_{H^2}+\|\tilde{Z}_j\|_{H^1} > \epsilon$$
    for any affine change $(\tilde{H}_j,\tilde{Z}_j)$. Fix a convergent subsequence of $\tau_j$ or a subsequence $\tau_j\to -\infty$. From Lemma \ref{Z_error},
$$Z_j = a_j\cdot n(\theta)+o_{H^1}(1).$$
Also note that $W_j = H_j+Z_j>0$. Since $H_j$'s satisfy $c_1 \leq H_j \leq C_1$, $a_j$'s have to be bounded. If $\inf W_j = m_j>0$ then since $H^1$ bound gives $C^{0,1/2}$ bound, using the fact that $W_j \leq m_j+M|\theta-\theta_0|^{1/2}$ , we get $-\log(m_j)\leq CI_j$. Since the monotone quantity $Q$ was bounded, by \eqref{MQ} we get that $m_j\geq c> 0  $ for some constant $c$. If we set
$$
c_j:=\frac{2\pi}{I_j}.
$$
Then by $E_2\to 0$,
$$
\|R_jW_j^2-c_jW_j^{-1}\|_{L^2}\to0.
$$
By boundedness of $W_j$, $R_j$ is uniformly bounded in $L^2$. Since $R_j=H_j+H_{j,\theta\theta}$ and $H_j$ is bounded in $L^2$, $H_j$'s are uniformly bounded in $H^2$. Passing to a subsequence,
$$
H_j\to H_\infty\quad\text{strongly in }H^1\text{ and }C^0,
\qquad
a_j\to a,
$$
and
$$
W_j\to W_\infty:=H_\infty+a\cdot n.
$$
Also $c_j\to c>0$. Passing to the limit in the $E_2$ identity gives
\begin{equation}
R_\infty W_\infty^3=c.
\end{equation}
Adding the $k=1$ function $a\cdot n$ does not change $R=H+H''$. Therefore, if
$$
H_e:=H_\infty+a\cdot n=W_\infty,
$$
then
\begin{equation}
H_{e,\theta\theta}+H_e=cH_e^{-3}.
\end{equation}
From classical affine geometry (see, for example, \cite{A1}), the only solutions to the above ODE are support functions of ellipses. Therefore up to an affine change $L$, $H_{\infty} =1$ and $Z_{\infty} = 0$. Also note that the entries of $L$ are completely governed by the fact that $B_{c_1} \subseteq E_{\infty}\subseteq B_{C_1}$ for the limiting ellipse $E_{\infty}$ and the translation term $a\cdot n$. Thus there exists an affine change such that $(\tilde{H}_j,\tilde{Z_j})$ satisfy \eqref{close_circle}, a contradiction. 
\end{proof}

\begin{rmk}\label{john_rmk}
    Note that the above proof shows that for John normalized sections $K_{\tau}$, since the total dissipation is finite by \eqref{finite_dissipation}, there exists a subsequence $\tau_j$ such that after another fixed affine change $L \in ASL(2)$, the sections $LK_{\tau_j}$ converge to the unit circle.
\end{rmk}

\begin{rmk}\label{radial_rigidity}
    From the proof of the preceding lemma it follows immediately that if $\delta_1 = \delta_2 = 0$ then the section at height $t$ is an ellipse.
\end{rmk}

Now we have to establish 2 things. Stability, i.e., once the normalized sections are near a circle, their distance from the circle is controlled by the dissipations. Another is that we need to bound the $SL(2)$ components of the $ASL(2)$ deformations we make to ensure uniform eccentricity of sections.

\subsection{Dissipations Control Distance from Circles} We establish that if an area normalized section $K_{\tau}$ is near the circle in some strong sense then the dissipation is bounded below by the  distance of the convex body $K_{\tau}$ from the circle in the strong sense.

\begin{lemma}\label{strong_small}
    Let
    $$
H=1+\phi.
$$
Further assume the strong smallness condition
\begin{equation}
X(\tau):=\|\phi\|_{H^2}^2+\|Z\|_{H^1}^2\le\varepsilon_{ss}
\end{equation}
for some $\varepsilon_{ss}$ small. If $P_1\phi = P_2\phi = 0$ then
$$X(\tau) \leq C\mathcal{E}(\tau).$$ 

\end{lemma}
\begin{proof}  Expand the constant area condition \eqref{fix_area},
$$
2\ic\phi d\theta+\ic(\phi^2-\phi_\theta^2)d\theta=0.
$$
Therefore
\begin{equation}
\left|\ic\phi d\theta\right|\leq C\|\phi\|_{H^1}^2.
\end{equation}
Similarly \eqref{orthogonality} gives
$$
\ic Z=-\ic(\phi Z-\phi_\theta Z_\theta).
$$
So
\begin{equation}\label{z_0_small}
\left|\ic Z\right|\leq C\|\phi\|_{H^1}\|Z\|_{H^1}.
\end{equation}
By the definition of $E_1$,
$$E_1 = \sum_{k=2}^{\infty}(k^2-1)|Z_k|^2-|Z_0|^2.$$
Therefore by \eqref{z_0_small},
\begin{equation}
\|P_{\geq2}Z\|_{H^1}^2\leq CE_1+CX^2.
\end{equation}
Let
$$
F:=RW^2-\frac{2\pi}{I}W^{-1}.
$$
Then $E_2=\|F\|_{L^2}^2$. Put $s:=\phi+Z$, so $W=1+s$. By Sobolev embedding in one dimension, $\|s\|_{C^{0,1/2}}\leq C\sqrt{\varepsilon_{ss}}$. The preceding estimates also imply
$$
\ic s=O(X).
$$
Therefore
$$
I=\int(1+s)^{-2}=2\pi+O(X),
\qquad
\frac{2\pi}{I}=1+O(X).
$$
Now
$$
\begin{aligned}
RW^2
&=(1+\phi+\phi_{\theta\theta})(1+s)^2\\
&=1+(\phi+\phi_{\theta\theta})+2s+O_{L^2}(X)\\
&=1+(3\phi+\phi_{\theta\theta})+2Z+O_{L^2}(X),
\end{aligned}
$$
while
$$
\frac{2\pi}{I}W^{-1}=1-\phi-Z+O_{H^1}(X^2).
$$
Consequently
\begin{equation}\label{F_equation}
F=(\partial_{\theta\theta}+4)\phi+3Z+\mathcal{R},
\qquad
\|\mathcal{R}\|_{L^2}\leq CX.
\end{equation}

Projecting equation \eqref{F_equation} to $k=1$ mode we get,
$$
3P_1Z=P_1F+O_{L^2}(X),
$$
since $P_1\phi=0$. Therefore,
\begin{equation}
\|P_1Z\|_{H^1}\le CE_2^{1/2}+CX.
\end{equation}

The operator $L_0:=\partial_{\theta\theta}+4$ has kernel exactly the $k=2$ space. Since $P_2\phi=0$,
$$
\|\phi\|_{H^2}\leq C\|L_0\phi\|_{L^2}.
$$
Therefore we get,
\begin{equation}
\|\phi\|_{H^2}\leq C\mathcal{E}^{1/2}+CX.
\end{equation}

Thus,
$$\sqrt{X} \leq C\mathcal{E}^{1/2}+CX.$$
For $\varepsilon$ small we get,
$$X \leq C\mathcal{E}.$$

\end{proof}

\subsection{Gauge Fixing} 
Since we can only get the bounds in Lemma \ref{strong_small} under the conditions $P_1H=P_2H=0$, we have to gauge fix the sections $K_\tau$. While $P_1H$ can be made zero by simple translations of $K_\tau$, we now show that for $\|\phi\|_{C^0}$ small, we can in fact find $SL(2)$ deformations near the identity to ensure $P_2H = 0$.\\
Recall that if $h$ is the support function corresponding to a convex body $K$, then $\tilde{h}=h+a\cdot n(\theta)$ is the support function for $\tilde{K}=K+a$. Given $S\in SL(2)$, we consider $\tilde{K} = SK$. Then $\tilde{h}(x) = h(S^{T}x)$ is the support function for $\tilde{K}$. Here we have taken a 1-homogeneous extension of both the support functions. In particular near $I_{2\times 2}$, we write $A = \exp(\e M)$ for some $M\in\mathfrak{sl}(2)$. At fixed normal $n(\theta)$, the first variation of the support function is
\begin{equation}
\xi_Mh=n\cdot Mx,\qquad x=hn+h_\theta n^{\perp}.
\end{equation}
Thus
$$
\xi_Mh=a_M(\theta)h+b_M(\theta)h_\theta,
$$
where
$$
a_M=n\cdot Mn,\qquad b_M=n\cdot Mn^{\perp}.
$$
\begin{lemma}\label{IFT}
    Let $H = 1+\phi$ be the support function of a convex body $K$. There exists $\varepsilon_{GF}$ such that for $\|\phi\|_{C^0} \leq \varepsilon_{GF}$ there exists $L(H) \in ASL(2)$ such that the support function of the convex body $L(H)K$, $\tilde{H}$, satisfies,
    $$P_1\tilde{H} = P_2\tilde{H} = 0$$
    and the mapping $H \mapsto L(H)$ is $C^1$ in the ball $\|\phi\|_{C^0} \leq \varepsilon_{GF}$. In fact $L(H) = (\exp(M(H)),c(H))$ for some $M(H)\in \mathfrak{sl}(2)$ and $c(H) \in \mathbb{R}^2$.
\end{lemma}

\begin{proof}
    Write
$$
c=(c_1,c_2),
\qquad
M=M(a,b)
=
\begin{pmatrix}
a&b\\
b&-a
\end{pmatrix}.
$$
Identify the first two Fourier modes with $\mathbb{R}^4$ by
$$
\mathcal{C}(H)
:=
\begin{pmatrix}
\frac1\pi\int_0^{2\pi}H(\theta)\cos\theta\,d\theta\\
\frac1\pi\int_0^{2\pi}H(\theta)\sin\theta\,d\theta\\
\frac1\pi\int_0^{2\pi}H(\theta)\cos2\theta\,d\theta\\
\frac1\pi\int_0^{2\pi}H(\theta)\sin2\theta\,d\theta
\end{pmatrix}.
$$
Then
$$
\mathcal C(H)=0
\quad\Longleftrightarrow\quad
P_1H=P_2H=0.
$$
Define
$$
\mathcal{F}(c_1,c_2,a,b;H)
:=
\mathcal C(T_{c,M(a,b)}H),
$$
where
\begin{equation}
T_{c,M}H(\theta)
=
|e^M\theta|H\left(\frac{e^{M}\theta}{|e^{M}\theta|}\right)
+
c\cdot n(\theta).
\label{eq:affine-action}
\end{equation}
We want to solve
$$
\mathcal{F}(c_1,c_2,a,b;H)=0.
$$
At the disk,
$$
\mathcal{F}(0,0,0,0;1)=0.
$$
While the mapping $M\mapsto T_{c,M}(H)$ need not be differentiable, however the Fourier projections can be written linearly in terms of $H$. For some trigonometric polynomial $q$,
$$
\int_0^{2\pi}
|e^M\theta|H\left(\frac{e^{M}\theta}{|e^{M}\theta|}\right)q(\theta)\,d\theta
=
\int_0^{2\pi}
H(\phi)K_q(M,\phi)\,d\phi.
$$

Here $K_q$ is a smooth function of $(M,\phi)$ near $M = 0$. Therefore $\mathcal{F}$ is $C^1$ and in fact linear in $H$. Computing the partial derivative $D_{c=(0,0);M = (0,0)}\mathcal{F}$ at the function $1$ we get,
$$
D_c\mathcal{F}((0,0),(0,0);1)[c_1,c_2]
=
\begin{pmatrix}
c_1\\
c_2\\
0\\
0
\end{pmatrix}, \;\;
D_{(a,b)}\mathcal{F}((0,0),(0,0);1)[a,b]
=
\begin{pmatrix}
0\\
0\\
a\\
b
\end{pmatrix}.
$$
Thus $D_{c,M}\mathcal{F}((0,0),(0,0);1) = Id_{\mathbb{R}^4}$. Therefore by implicit function theorem we get $C^1$ mapping and $\varepsilon_{GF}$ small.

\end{proof}

We have the following important consequence of the above lemma,

\begin{lemma}\label{drift_control}
    Let $K_t$ be convex bodies evolving by \eqref{2nd_order_evol}. Let $(H(\tau,\cdot),Z(\tau,\cdot))$ be the support function and the physical velocity fixed area convex bodies $K_{\tau}$ corresponding to $K_t$. If $\|H(\tau_0)-1\|_{C^0} \leq \varepsilon_{1}$ then for some $\tau_1<\tau_0$ there exists a $C^1$ mapping $\tau\mapsto L(\tau)\in ASL(2)$ on the interval $(\tau_1,\tau_0)$. The support function, $\tilde{H}(\tau)$ for the affine normalized sections $L(\tau)K_{\tau}$ has the following properties,
    $$P_1\tilde{H} = P_2\tilde{H} = 0.$$
    Moreover for the $SL(2)$ part of $L(\tau)$, $M(\tau)$, we have,
    $$|M'(\tau)(M(\tau))^{-1}| \leq C|P_2\tilde{Z}|.$$
    Here $\tilde{Z}$ is the physical velocity for the section $L(\tau)K_{\tau}$. 
    
\end{lemma}
\begin{proof}
    By Lemma \ref{IFT}, for $\varepsilon_1$ small, we have maps $L(\tau) = L(H(\tau)) = (M(\tau),c(\tau))$ for some time $(\tau_1,\tau_0)$ which are $C^1$ since the convex bodies $K_{\tau}$ evolve in a $C^{1,1}$ fashion. Therefore for the affine normalized convex bodies we have $P_1\tilde{H} = P_2\tilde{H} = 0$. For the quantities $(\tilde{H},\tilde{Z})$ corresponding to $L(\tau)K_{\tau}$ note,
    $$\frac{d}{d\tau}\tilde{H} = \tilde{Z}+\xi_\Omega \tilde{H}+(k=1\text{ translation mode})$$
    for $\Omega = M'M(\tau)^{-1}$. Projecting onto $k=1$ Fourier mode,
$$0 = P_1(\tilde{Z}+\xi_{\Omega}\tilde{H})+(k=1 \text{ translation mode})$$
Therefore,
\begin{equation}\label{total_derivative}
    \frac{d}{d\tau}\tilde{H} = (I-P_1)(\tilde{Z}+\xi_{\Omega}\tilde{H})
\end{equation}
Projecting onto $k=2$ Fourier mode,

\begin{equation}
0=P_2\tilde{Z}+P_2(\xi_\Omega \tilde{H}).
\end{equation}
At the disk $\tilde{H}\equiv1$, for
$$
\Omega=
\begin{pmatrix}
a&b\\ b&-a
\end{pmatrix},
$$
we have
$$
\xi_\Omega1=n\cdot\Omega n=a\cos2\theta+b\sin2\theta.
$$
We note that change due to the skew-symmetric part of $\Omega$ for a $M$ can be controlled by the symmetric part at the disc. Thus
$$
\Omega\longmapsto P_2(\xi_\Omega1)
$$
is an isomorphism from symmetric trace-free matrices to the $k=2$ Fourier space. For a general h we have,

$$ P_2(\xi_\Omega h) = T_1(a,b,P_0h,P_4h)\cos(2\theta)+T_2(a,b,P_0h,P_4h)\sin(2\theta)$$
where $T_1,T_2$ are linear in $(a,b)$ as well as $h$. By perturbation it remains uniformly invertible when $H$ is sufficiently close to the disk. Consequently
\begin{equation}\label{guage_fix}
|\Omega|\le C|P_2\tilde{Z}|
\end{equation}
whenever $\|\tilde{H}-1\|_{C^{0}} \leq \varepsilon_{2}$. Choose $\varepsilon_1$ small such that it's less than $\varepsilon_{GF}$ and such that the resulting $\tilde{H}$ are such that $\|\tilde{H}-1\|_{C^0} \leq \|\tilde{H}-H\|_{C^0}+\|H-1\|_{C^0} \leq \varepsilon_2$.
\end{proof}

\begin{rmk}\label{translation_small}
    Note that above proof also gives us bounds on the $k=1$ translation term in terms of $P_1Z$ and $P_1H$ after gauge fixing.
\end{rmk}

\subsection{Affine Drift Control}
 If we are in the setting of Lemma \ref{drift_control} we have some control on the affine drift $M'$. However note that apriori we only have $|P_2Z| \leq \sqrt{\mathcal{E}}$ and we don't know its integrability. We would therefore like to upgrade the bound to $\mathcal{E}$. Classically for geometric flows such improvements in exponents are usually done by Łojasiewicz inequalities for e.g. in \cite{Simon}, \cite{CM1}. However since our flow also has a nontrivial variational structure \eqref{EL} with affine symmetry, we can exploit that to get better bound. Consider the Noether momentum corresponding to the functional \eqref{EL} and the infinitesimal symmetry $\exp(\epsilon M)$ for $M \in \mathfrak{sl}(2)$,
$$J(t) = \ic v^{-3}\xi_M h d\theta.$$

We can write,
$$J'(t) = \ic v^{-3}\xi_M v d\theta+\ic (v^{-3})_t\xi_M h d\theta = -\frac{1}{2}\frac{d}{d\epsilon}G(\epsilon) + 3f\frac{d}{d\epsilon}A(\epsilon).$$
We used \eqref{2nd_order_evol} for the second term. Here $G(\e), A(\e)$ are quantities associated to the convex bodies $\exp(\e M)K_t = K_{\e,t}$. Since we know they are affine invariant we get 
$$J'(t) = 0.$$
Therefore, $J_M(t) = J_M(0) = 0$. Now using our area normalized quantities we have,
$$\ic W^{-3}\xi_M H d\theta = 0 \qquad \forall M \in \mathfrak{sl}(2).$$
We use the above equality to control the $k=2$ harmonics of $Z$.
\begin{lemma}\label{good_drift_bound}
    Assume the strong smallness condition
\begin{equation}
X:=\|\phi\|_{H^2}^2+\|Z\|_{H^1}^2\leq \varepsilon_{ss}
\end{equation}
for $\varepsilon_{ss}$ small as in Lemma \ref{strong_small}. If $P_1\phi = P_2\phi = 0$ then
$$|P_2Z| \leq C \mathcal{E}$$
\end{lemma}

\begin{proof}
    
Take
$$
S_1=
\begin{pmatrix}
1&0\\0&-1
\end{pmatrix}.
$$
Then
$$
\xi_{S_1}H=H\cos2\theta-H_\theta\sin2\theta,
\qquad
\xi_{S_1}1=\cos2\theta.
$$
The zero-momentum identity is
$$
0=\ic W^{-3}\xi_{S_1}H\,d\theta.
$$
Using $W = 1+\phi+Z$,
$$
W^{-3}=1-3(\phi+Z)+O_{C^0}(X).
$$
Also
$$
\xi_{S_1}H=\cos2\theta+\xi_{S_1}\phi.
$$
Note that,
$$
\ic\xi_{S_1}\phi\,d\theta
=\ic(\phi\cos2\theta-\phi_\theta\sin2\theta)\,d\theta
=3\ic\phi\cos2\theta\,d\theta.
$$
Hence
$$
\ic\xi_{S_1}\phi-3\ic\phi\cos2\theta=0.
$$
Therefore,
\begin{equation}
    \begin{aligned}
        0 = \ic W^{-3}\xi_{S_1}H\,&d\theta =\ic\xi_{S_1}\phi-3\ic\phi\cos2\theta\\
        & -3\ic Z\cos2\theta+O(X).
    \end{aligned}
\end{equation}
 Thus
$$
\left|\ic Z\cos2\theta\,d\theta\right|\leq CX.
$$
For
$$
S_2=
\begin{pmatrix}
0&1\\1&0
\end{pmatrix}
$$
one similarly gets the sine component. Therefore
\begin{equation}
|P_2Z|\leq CX\leq C\mathcal{E}.
\end{equation}

\end{proof}

\subsection{Stability Near Circles}
We have now established being $C^0$ near a circle is enough to gauge fix. Moreover being close to the circle in the strong sense gives us lower bounds for the dissipations as well as ensures good bounds for the affine drifts. We now use a standard  argument to show we can never exit the strong small neighbourhood as $\tau\to -\infty$.\\

 Let us John normalize every section $K_{\tau}$ and consider $(H,Z)$ for the normalized sections. Let $\tau_0$ be a time when $\mathcal{E} \leq \delta_0$ such that by Lemma \ref{normalization}, for the affine normalized section $L(\tau_0)K_{\tau_0}$, $X(\tau_0) \leq \varepsilon_*/4$ for $\varepsilon_* \leq \varepsilon_{ss}$ as in Lemma \ref{strong_small}. Further normalize every section after this with $L = L(\tau_0)$ as in Lemma \ref{normalization}. Let $(\tau_1,\tau_0]$ be the maximal semi-open interval in which we could ensure $P_1H=P_2H=0$ by $ASL(2)$ transformations. After gauge fixing consider the distance function,
$$d(\tau) = \|H(\tau)-1\|_{C^0}.$$
 Let $(\tau_*,\tau_0]\subseteq (\tau_1,\tau_0]$ be the maximal interval such that  $d < \varepsilon_*/2$. Note that by Lemma \ref{IFT} we must have $\tau_*>\tau_1$ if $\tau_*>-\infty$. Assume not.  For $\varepsilon_*$ small we ensure $\|H(\tau)-1\| \leq \varepsilon_{GF}/2$ for all $\tau \in (\tau_*,\tau_0]$. Therefore the gauges $L(\tau)$ are bounded uniformly in $SL(2)$ and have a subsequence, $L(\tau_j)$, converging to some $L^*$. Since $K_\tau$ were evolving continuously, along the subsequence, $L(\tau_j)K_{\tau_j}$, we have $L(\tau_j)K_{\tau_j}\to L^*K_{\tau_*}$ such that the support functions have $P_1H(\tau_j)=P_2H(\tau_j) = 0$. Hence at $\tau = \tau_*$ we can still gauge fix $K_{\tau_*}$ using an $ASL(2)$ transformation, a contradiction.  
 We also note a compactness result relating $d(\tau)$ and the dissipation $\mathcal{E}$,

 \begin{lemma}\label{small_C^0}
    For every $\epsilon>0$, there exists $\delta>0$ such that if $\mathcal{E}(\tau) \leq \delta$ ,$d(\tau) \leq \varepsilon_*$ and $P_1H=P_2H=0$ then
    $$\|H-1\|_{H^2}+\|Z\|_{H^1} \leq \epsilon$$
\end{lemma}
\begin{proof}
    $d(\tau) \leq \varepsilon_*$ implies $B_{1-C\varepsilon_{*}} \subseteq K_\tau \subseteq B_{1+C\varepsilon_{*}}$. By Lemma \ref{normalization} we have the strong smallness conclusion under some $ASL(2)$ transformation. However we also note that $P_1H = P_2H=0$. Moreover note that in the proof of Lemma \ref{normalization}, the support function of the limiting ellipse $E_{\infty}$ must satisfy $P_1H_{\infty} = P_2H_{\infty}=0$. This is only possible if $E_{\infty}$ is the unit circle. Therefore the $ASL(2)$ deformation must be identity map.
\end{proof}

We choose $\epsilon = \varepsilon_{ss}$ as in Lemma \ref{strong_small} to get $\delta_*$ from above Lemma. We now restrict ourselves to the interval $(\tau_*,\tau_0]$. If $\mathcal{E} \leq \delta_*$ then,

\begin{equation}\label{c_bound}
    |d(\tau)|^2 \leq X \leq C\mathcal{E}.
\end{equation}
Else the bound is trivially true for C $\geq \varepsilon_*/\delta_*$. Moreover since the map $\tau\mapsto H(\tau)$ is $C^1$, we get by \eqref{total_derivative},
$$|d'(\tau)| \leq \|H_\tau\|_{C^0} \leq \|(I-P_1)(Z+\xi_{\Omega}H)\|_{C^0}.$$

Since $d(\tau) \leq \varepsilon_*$, we have $B_{1-C\varepsilon_{*}} \subseteq K_\tau \subseteq B_{1+C\varepsilon_{*}}$. By Lemma \ref{Z_error} we have $\|(I-P_1)Z\|_{H^1}^2 \leq C E_1$. Also note
 $$|\xi_{\Omega}H| \leq |\Omega||x(\theta)| \leq C|\Omega| \leq C|P_2Z|.$$
 Therefore $\|(I-P_1)\xi_{\Omega}H\|_{C^0} \leq C E_1^{\frac{1}{2}}$. Thus we have,
 $$|d'(\tau)| \leq C\mathcal{E}^{1/2}.$$
 \begin{lemma}\label{exit_cost}
If $\tau_*>-\infty$ then there is $\eta_*>0$ such that on the interval $(\tau_*,\tau_0]$
$$
\int \mathcal{E}(\tau)\,d\tau\ge\eta_*.
$$
\end{lemma}

\begin{proof}
Consider a subinterval $J = (\tau_*,\tau_1] \subseteq (\tau_*,\tau_0]$ near $\tau_*$ such that the following bounds hold
$$
\frac{\varepsilon_*}{4}\le d\le\frac{\varepsilon_*}{2}.
$$
Moreover we also assume $d(\tau_1) = \epsilon_*/4$. By \eqref{c_bound},
$$
\mathcal{E}\geq c_0\varepsilon_*^2\qquad\text{on }\;J.
$$
If $A_*:=\int_{J}\mathcal{E}$, then
$$
|J|\le\frac{A_*}{c_0\varepsilon_*^2}.
$$

Moreover,
$$
\frac{\varepsilon_*}{4}
\leq\int_{J}|d'|\,d\tau
\leq C\int_{J}\mathcal{E}^{1/2}\,d\tau
\leq C|J|^{1/2}A_*^{1/2}
\leq\frac{C}{\varepsilon_*}A_*.
$$
Therefore $A_*\geq \tilde{c}\varepsilon_*^2=: \eta_*>0$.
\end{proof}

\begin{proof}[Proof of Theorem \ref{x=0_theorem}]
    
Choose a $\tau_0$ small such that $\mathcal{E}(\tau_0) \leq \delta_0$. This $\delta_0$ is chosen such that after affine normalization, $L(\tau_0)$ according to Lemma \ref{normalization}, we get $X(\tau_0) \leq \mu$ for some $\mu$ small. We first note that $Q(-\infty) = \lim_{\tau\downarrow -\infty}Q(\tau) = 8\pi^2$. By Remark \ref{john_rmk}, we know along a subsequence our sections after affine normalization converge to $(H,Z) = (1,0)$. Therefore $W=1$ and we get $I = \ic W^{-2}d\theta = 2\pi$. Therefore $Q(-\infty) = 8\pi^2$. Now also note that,
$$\int_{-\infty}^{\tau} \mathcal{E}(s)ds \leq  C\log\left(\frac{Q(-\infty)}{Q(\tau)}\right).$$
Since we have $X(\tau_0) \leq \mu$. We choose $\mu$ small such that,
$$\log\left(\frac{Q(-\infty)}{Q(\tau_0)}\right) \leq \frac{\delta_*}{C}.$$
 The choice of $\mu$ and hence $\delta_0$ were made so that both $\mathcal{E}(\tau_0)$ and $\int_{-\infty}^{\tau_0}\mathcal{E}d\tau$ are less than $ \eta = \min\{\eta_*/2,\delta_*\}$. Therefore as discussed above $X(\tau_0) \leq \varepsilon_*/4$ for such a choice. Note that since $\int \mathcal{E}\leq M$, for $M$ universal. Therefore we have $|\{E\geq \delta_0\}| \leq M/\delta_0$. Therefore we can choose a $\tau_0\geq -M/\delta_0$ universally away from $-\infty$. Normalize all sections $\tau \leq \tau_0$ by $L = L(\tau_0)$. Therefore by Lemma \ref{exit_cost}, $\tau_1 = \tau_* = -\infty$. On the interval $(-\infty,\tau_0]$,
$$|\Omega| = |M'M^{-1}| \leq C|P_2Z| \leq C\mathcal{E}$$
by dividing into two cases $\mathcal{E}\leq \eta$ and $\mathcal{E} \geq \eta$. If $\mathcal{E} \leq \eta$ then use Lemma \ref{small_C^0} along with Lemma \ref{good_drift_bound}. If $\mathcal{E} \geq \eta$ use the trivial bound as in Lemma \ref{Z_error}. Now using integrability of the dissipations,
\begin{equation}
\int_{-\infty}^{\tau_0}|\Omega(\tau)|\,d\tau<\infty.
\end{equation}
Since $M'M^{-1}=\Omega$, then for $a<b\leq\tau_0$,
$$
\|M(b)M(a)^{-1}-I\|
\le\exp\left(\int_{a}^{b}|\Omega|\right)-1.
$$
Hence the function $\tau \mapsto M(\tau)$ is bounded and uniformly continuous. Therefore $M(\tau)\to M_{-\infty}$ as $\tau\to-\infty$. Therefore since $B_{1-C\varepsilon_*} \subseteq L(\tau)K_{\tau} \subseteq B_{1+C\varepsilon_*}$, we get that all $K_\tau$ have bounded eccentricity and in fact $K_\tau$ converge to an ellipse. 

\end{proof}

\begin{rmk}
    From the above proof as well as Remark \ref{translation_small} we can also control translation terms. Since the integral of dissipations go to $0$, translations converge to $0$.
\end{rmk}

\section{Applications}\label{application}

\subsection{Affine rigidity} As an application of the above technique we establish the following rigidity of radial solutions,

\begin{theorem}\label{rigidity}
Let $\Omega=B_R(0)\subset\mathbb{R}^2$, and let $u\in C^2(\overline{\Omega}\setminus\{0\})$ be a convex solution of
$$
\det D^2u=f(u)
\quad\hbox{in }\Omega,
\qquad
u=1\quad\hbox{on }\partial\Omega,
$$
with an isolated minimum at $x=0$ with $u(0)=0$. Assume
$$
c_ft^p\le f(t)\le C_ft^p
$$
for some $-2< p<2$ and $t\leq 1$.  Then every section is a concentric disk, and
$$
u(x)=U(|x|)
$$
for some increasing radial profile $U$.
\end{theorem}
\begin{proof}
    By Lemma \ref{monotone}, the quantity $Q$ is monotone and by Lemma \ref{aff_inv} it's affine invariant. By Theorem \ref{x=0_theorem}, $r(t)L(t)S_t \to B_1(0)$. Note that we know from Section \ref{rescale} that for
    $$I = \ic W^{-2}d\theta,$$
    $Q = 4\pi I$ by \eqref{MQ}. We also know that after an affine change $W(\tau) \to 1$. Therefore,
    $$Q(t) \leq 8\pi^2$$
    for all $t>0$. At $t=1$ we have
    $$Q(1) = \frac{1}{A}\left(\ic v^{-2}d\theta\right) \left(\ic \rho vd\theta\right)^{2} = \frac{1}{\pi }\left(\ic v^{-2}d\theta\right) \left(\ic  vd\theta\right)^{2},$$
    for $v = h_t$. Here $h(t,\cdot)$ is the support function corresponding to $S_t$. By Jensen's inequality,
    $$\left(\ic v\frac{d\theta}{2\pi}\right)^{-2} \leq \ic v^{-2}\frac{d\theta}{2\pi}.$$
    Therefore $Q(1) \geq \frac{1}{\pi}8\pi^3 = 8\pi^2$. Thus $Q(t)$ has to be constant and by Lemma \ref{monotone} the defects $\delta_1,\delta_2$ have to be identically $0$. The conclusion follows from the equality case discussed in Lemma \ref{normalization}.
\end{proof}

\begin{rmk}
    Such a rigidity result was obtained in Cheng-Dai-Li \cite{DCL} under the assumption that the right hand side is differentiable using the moving plane method. On the plane, the moving plane method was used in Liu-Ma \cite{LM} to give symmetry results for a general class of right hand side under monotonicity assumptions of the right hand side.
\end{rmk}

Note that while we analyzed sections $S_t$ as $t\to 0$, one could do the same as $t\to \infty$. We now prove Theorem \ref{Liouville2},
\begin{proof}[Proof of Theorem \ref{Liouville2}]
    Note that Lemma \ref{rescaling}, the bounds on $Q$ and the dissipation relation \eqref{dissipation_rel} remain finite. Therefore the sections $S_t$ as $t\to \infty$ also have the property $r(t)L(t)S(t)\to B_1(0)$ as in Theorem \ref{x=0_theorem}. Proof follows exactly as the proof of Theorem \ref{x=0_theorem}. Therefore as in the proof of Theorem \ref{rigidity}, not only $Q(0^{+}) = 8\pi^2$, we must also have $\lim_{t\to \infty}Q(t) = 8\pi^2$. By monotonicity of $Q$, we get $\delta_1=\delta_2=0$ for all $t$. Therefore radial rigidity follows from the proof of Lemma \ref{normalization}.
\end{proof}

\subsection{Classification Results for a Free Boundary Problem}
For the special case $f_q(t) = t^q$ and $0\leq q<2$, we can classify all the finite measure coincidence sets of the free boundary problem 
\begin{equation}\label{fbp}
    \det(D^2u) = u^q\chi_{\{u>0\}} \text{ on } \mathbb{R}^2.
\end{equation}
 We can also classify all global solutions to the corresponding equation. This completes the program started by Jin-Tu-Xiong in \cite{JTX} to classify all global solutions the the free boundary problem \eqref{fbp} in $2$ dimensions. Note that the free boundary is denoted by $K = \{u=0\}$. We recall an important result from \cite[Theorem 1.5]{JTX}.

\begin{theorem}
Let $n\ge 2$ and $0\le q<n$. Suppose $v$ is a solution to
\begin{equation}\label{eq:obstacle problem}
	\det D^2 v = g(x,v,Dv) v^q\chi_{\{v>0\}},
\end{equation}
and $0$ is a strictly convex point on the free boundary. Assume the free boundary can be locally written as the graph of the function $\psi$ near $0$. Further assume that $ 0<\lambda \leq g \leq \Lambda <\infty$, and $\left\|g\right\|_{C^{\alpha}} \leq \Lambda$ for $\alpha \in (0,1)$. Then, we have
\[
\left\|\psi\right\|_{C^{2,\beta} (B_c'(0))} \leq C  ,\quad 0 < c I_{n-1}\leq D_{x'}^2 \psi \leq  C  I_{n-1}\quad  \text{in }B_c'(0),
\]
and
\begin{equation}\label{eq:v distance optimal}
c\operatorname{dist}(x,K)^{\frac{n+1}{n-q}} \leq v(x) \leq C\operatorname{dist}(x,K)^{\frac{n+1}{n-q}}   \quad \text{in }B_{c}(0) ,
\end{equation}
where $\beta=\min\left\{\frac{2(q+1)}{n-q},1\right\}\alpha$,  $c$ and $C$ are positive constants depending only on $n$, $q$, $\alpha$, $\lambda$ and $\Lambda$.
Moreover, if $g(x,v,Dv)=g(x,v)$ is independent of $Dv$, then we can take $\beta=\alpha$.   
\end{theorem}

We also recall that for a global solution \eqref{fbp}, the entire free boundary has to be strictly convex by \cite[Proposition 2.18]{JTX}. Then from \cite[Section 5]{JTX} the free boundary is $C^{2,\alpha}$. By performing a partial Legendre transform, $u^*(p,y)$, in the $x$ variable we define,
 $$\psi = \frac{u^*}{p}.$$
As $p\to 0^{-}$ we get, $\psi(0,y) = \partial_p u^*(0^{-},y) = x$. Therefore $(y,\psi(0,y))$ is the graphical representation of the free boundary locally. Moreover using $u^* = p\psi$ we get,
\begin{equation}\label{u_psi_p}
    u = p u^*_p - u^* = p^2\psi_p.
\end{equation}
We restrict ourselves to a subset of the half plane $\mathbb{H} = \{p<0\}$. By a change of variable, $p = -s^{2/k}$ for $k = \frac{n-q}{q+1}$, the PLT equation gives us,
    \begin{equation}\label{psi_plt_eqn}
        \psi_{ss}+\frac{k+2}{k}\frac{\psi_s}{s}+\left(\frac{2}{k}\right)^{q+2}
\left(-\frac{\psi_s}{s}\right)^{-q}
\psi_{yy}
=0 \;\text{ on } B_c^+(0).
    \end{equation}
  Here $B_c^+(0) = B_c(0)\cap \{s>0\}$ for a small radius $c$. We now recall 2 important results pertaining to $\psi$ from \cite[Lemma 5.1, Theorem 5.3]{JTX}.

\begin{lemma}\label{unif_control_psi}
 The function $\psi$ satisfies 
\begin{equation}
c\leq 	-\frac{\psi_s}{s}  \leq C  \quad  \text{and }\quad 	c \leq   \psi_{yy} \leq C  \quad  \text{in }B_c^+(0).
\end{equation}
\end{lemma}

Note that Theorem $5.3$ in \cite{JTX} is stated slightly differently from below. In particular, they show only a part of the estimate below but for any $0<\alpha<1$. The problem arises because they can't control the $C^{\alpha}$ norm of $\frac{\psi_{sy}}{s}$ for any $\alpha$. In their proof, however, they note they do have control over the $C^{\alpha}$ norm of $\frac{\psi_{sy}}{s}$ for some small $\alpha_0$. Since we also require this estimate, we modify their theorem statement slightly.
\begin{theorem}\label{uniform_psi}
For some $\alpha \in (0,1)$, the function $\psi$ satisfies 
\begin{equation}
\left\|\psi  \right\|_{C^{2,\alpha} \left(\overline{B_c^+(0)}\right)} +  \left\| \frac{\psi_s}{s} \right\|_{C^{\alpha} \left(\overline{B_c^+(0)}\right)}+\left\| \frac{\psi_{sy}}{s} \right\|_{C^{\alpha} \left(\overline{B_c^+(0)}\right)} \leq C(n,q,\alpha)
\end{equation}
in the coordinates $(s,y)$.
\end{theorem}

We now use \eqref{u_psi_p} to control the level sets of $u$. Define,
$$a = -\frac{\psi_s}{s}.$$
By \eqref{u_psi_p} and $p\partial_p = \frac{k}{2}s\partial_s$, we get 
\begin{equation}\label{level_set_eqn}
    u = p^2 \psi_p = -\frac{k}{2}s^{1+\frac{2}{k}}\psi_s = \frac{k}{2}a s^m
\end{equation}
for $m = 2+2/k$. If we represent the level set $\{u=t\}$ locally as a graph in the $(s,y)$ variable as $(s_t(y),y)$, we get,
\begin{equation}\label{a_lvl_set_eqn}
    t = \frac{k}{2}a(s_t(y),y)(s_t(y))^m.
\end{equation}
By Lemma \ref{unif_control_psi}, $c \leq a \leq C$, therefore we have
\begin{equation}\label{uniform_s}
c\,t^{1/m}\le s_t(y)\le C\,t^{1/m}.
\end{equation}
Note that the actual local graphical representation of $\{u=t\}$ can be given by the relation,
$$x = \partial_p u^*.$$
Therefore writing everything in terms of $\psi$ we get,
$$x = \psi+p\psi_p = \psi+\frac{k}{2}s\psi_s.$$
Thus locally the graph of the level set $\{u=t\}$ is given by $(\gamma_t(y),y)$ where,
$$\gamma_t(y) = \psi+\frac{k}{2}s_t \psi_s.$$

We now calculate the derivatives of $\gamma_t$.

\begin{lemma}
    Along the level curve $\{u=t\}$,
\begin{equation}
\gamma_t'(y)=\psi_y(s_t(y),y) \;\text{ and }\; \gamma_t''(y) =\psi_{yy}(s_t,y)+\psi_{ys}(s_t,y)s_t'(y).
\end{equation}
\end{lemma}
\begin{proof}
    By \eqref{level_set_eqn},
    $$t = -\frac{k}{2}s_t^{1+2/k}\psi_s.$$
    Differentiating we get,
    $$0 = \left(1+\frac{2}{k}\right)s_t^{2/k}s_t'\psi_s+s_t^{1+\frac{2}{k}}\psi_s'.$$
    Therefore since
    $$\gamma' = \psi_y +s_t'\psi_s + \frac{k}{2}s_t'\psi_s + \frac{k}{2}s_t\psi_s'$$
    we immediately get $\gamma' = \psi_y$. The second equation follows from this.
    
\end{proof}

By \eqref{uniform_s}, we have $\gamma'_t\to \gamma_0'$ uniformly using Theorem \ref{uniform_psi}. Moreover in the expression of $\gamma_t''$ the $\psi_{yy}$ term also converges uniformly by Theorem \ref{uniform_psi}. We now control $s_t'$ using information about $a$. 
\begin{lemma}
The derivative $s_t'$ is bounded by $s_t$, i.e.,
    \begin{equation}
|\partial_{y}s_t(y)|\le C s_t(y).
\end{equation}
The constant is independent of $t$ for sufficiently small $t$.
\end{lemma}
\begin{proof}
    Differentiate \eqref{a_lvl_set_eqn},
    $$0=a_y s_t^m+ \left(a_s s_t^m+ma s_t^{m-1}\right)s_t'.
$$
Therefore,
$$s_t' = -\frac{a_ys_t}{a_ss_t+ma}.$$
By Theorem \ref{uniform_psi}, $a_y$ is bounded. Calculating $a_s$,
$$
a_s=-\frac{\psi_{ss}}{s}+\frac{\psi_s}{s^2}
=-\frac{\psi_{ss}}{s}-\frac{a}{s},
$$
so
$$
sa_s=-\psi_{ss}-a.
$$
Therefore
$$
ma+sa_s=(m-1)a-\psi_{ss}.
$$
Since $m-1=1+\frac2k=\frac{k+2}{k}$,
we have
\begin{equation}\label{denom}
ma+sa_s=\frac{k+2}{k}a-\psi_{ss}.
\end{equation}
Now rewrite \eqref{psi_plt_eqn} using $\psi_s/s=-a$,
$$
\psi_{ss}
-\frac{k+2}{k}a
+\left(\frac2k\right)^{q+2}a^{-q}\psi_{yy}
=0.
$$
Hence
$$
\frac{k+2}{k}a-\psi_{ss}
=
\left(\frac2k\right)^{q+2}a^{-q}\psi_{yy}.
$$
Therefore the denominator \eqref{denom} is bounded uniformly by Lemma \ref{unif_control_psi}.
\end{proof}

We now obtain the $C^2$ convergence of the level sets $\{u=t\}$ immediately,
\begin{corollary}
    For the solution to \eqref{fbp} the free boundary $\partial K$ is $C^{2,\alpha}$ and uniformly convex for some $\alpha>0$ when $0<|K|<\infty$. The nearby level curves
$$
\partial K_t=\{u=t\}
$$
converge to $\partial K$ in $C^2$ as $t\downarrow0$. In particular
\begin{equation}
\rho(t,\theta)=h(t,\theta)+h_{\theta\theta}(t,\theta)
\longrightarrow
\rho_0(\theta):=h_0+h_{0,\theta\theta}
\end{equation}
uniformly in $\theta$. Moreover,
\begin{equation}
0<c\le \rho_0(\theta)\le C<\infty.
\end{equation}
\end{corollary}

We now use our machinery to conclude the following result,

\begin{theorem}\label{classification}
Let $u\not\equiv 0$ be a convex solution of
$$
\det D^2u  = u^p\chi_{\{ u > 0\}},\quad u\geq 0,\quad  \text{in }\mathbb{R}^2.
$$
Let $K:=\{u=0\}$ denote the coincidence set. Assume $|K|<\infty$.  Then one of the following holds:
\begin{enumerate}
\item[(i).] If $|K| = 0$ then $K= \emptyset$ or is a singleton; 
\item[(ii).] If $0<|K|$ then $K$ is an ellipse. 
\end{enumerate} 
Moreover the solutions to \eqref{fbp} have to be affine equivalent to a radial solution in the above two cases.
\end{theorem}

The above result combined with \cite[Theorem 1.1]{JTX} classifies all global solutions and free boundaries to \eqref{fbp} in $2$ dimensions. If $|K|=0$ then \cite[Lemma 2.15, Proposition 2.18]{JTX} gives us $K$ has to be empty or singleton. Then using Theorem \ref{Liouville2} we conclude the affine equivalence to radial solutions. We now focus on the case $0<|K|<\infty$. To apply our monotonicity approach we need to ensure two things. First as $t\to \infty$, we need affine normalized sections to converge to the unit circle. Note that as long as we can ensure Lemma \ref{rescaling} holds in our setting we are done. This is indeed true by Remark \ref{rmk:large_t} since we are assuming $|K|<\infty$. Therefore by a similar proof to Theorem \ref{x=0_theorem} we get the convergence of affine normalized level sets to the unit circle. As argued in Theorem \ref{rigidity}, this also tells us $\lim_{t\uparrow \infty}Q(t) = 8\pi^2$. We therefore need to know $\lim_{t\downarrow 0}Q(t).$ We prove show the following result:
\begin{lemma}
    Let $u$ be a solution to \eqref{fbp} with $0<|K|<\infty$. Then
    $$\lim_{t\downarrow 0}Q(t) = \frac{L_{aff}(K)^3}{|K|}.$$
    Here $L_{aff}(K) = \int \kappa^{1/3}ds =  \ic \rho_0^{2/3}d\theta$ is the affine perimeter of the convex set $K$.
\end{lemma}
\begin{proof}
    By the $C^{1,\alpha}$ regularity of the solution to \eqref{fbp} we have,
    $$v^{-3} = |Du|^3 \to 0$$
    as $t\to 0$. Therefore by \eqref{2nd_order_evol} we get,
    $$v^{-3}(t,\theta) = 3\int_0^t s^q\rho(s,\theta)ds.$$
    Since we know the uniform convergence of $\rho(s,\theta)$ as $s\to 0$ we note,
    $$v^{-3}(t,\theta) = \frac{3}{q+1}t^{q+1}\rho_0(\theta)(1+o(1)).$$
    Therefore,
    $$v(t,\theta)=c_q t^{-\frac{q+1}{3}}\rho_0(\theta)^{-1/3}(1+o(1)),$$
    and
    $$v(t,\theta)^{-2} =c_q^{-2}t^{\frac{2(q+1)}3}\rho_0(\theta)^{2/3}(1+o(1))$$
    for $c_q=\left(\frac{q+1}{3}\right)^{1/3}$. Using these expansions first note that,
    $$G(t) = \ic v^{-2}d\theta=c_q^{-2}t^{\frac{2(q+1)}3}
\left(\ic\rho_0^{2/3}d\theta+o(1)\right).$$
Thus,
$$G(t)=c_q^{-2}L_{aff}(K)t^{\frac{2(q+1)}3}(1+o(1)).$$
We also get an expansion for $A'$ by,
$$A'(t)= \ic \rho vd\theta =c_q t^{-\frac{q+1}{3}}\left(\ic\rho_0^{2/3}\,d\theta+o(1)\right).
$$
Therefore,
$$A'(t)=c_qL_{aff}(K)t^{-\frac{q+1}{3}}(1+o(1)).$$
Combining the two calculations we get,
$$\lim_{t\downarrow 0}Q(t) = \lim_{t\downarrow 0}\frac{G(t)A'(t)^2}{A(t)} = \frac{L_{aff}^3(K)}{|K|}.$$
    
\end{proof}

\begin{proof}[Proof of Theorem \ref{classification}]
    By the classical affine isoperimetric inequality (see \cite{Sch}) 
    $$\lim_{t\downarrow 0} Q(t) = \frac{L_{aff}^3(K)}{|K|}  \leq 8\pi^2.$$
    Now using monotonicity of $Q$, $Q(t) \geq \lim_{t\to\infty}Q(t) = 8\pi^2$. Therefore by the arguments used in Theorem \ref{rigidity} and Theorem \ref{Liouville2}, $\delta_1 = \delta_2 = 0$ gives us radial solutions up to an affine transform and forces $K$ to an ellipse. 
\end{proof}

\begin{rmk}
    Note that in the proofs of Theorems \ref{rigidity}, \ref{Liouville2} and \ref{classification} we only required a sequence of level sets to converge to ellipses. This can be easily ensured by Remark \ref{john_rmk}. In particular, we did not require the full stability analysis carried out in Section \ref{sec:min_point_analysis}.
\end{rmk}

\section{Applications to Semigeostrophic Equations}\label{semigeostrophic}

The Semigeostrophic equations are given by,
\begin{equation}\label{SGE}
    \begin{split}
        \partial_t\nabla p_t+(\textbf{u}_t\cdot \nabla)\nabla p_t +\nabla^{\perp}p_t+\textbf{u}_t = 0&\;\text{ in }\; \mathbb{R}_{\geq 0}\times\Omega ,\\
        \text{div}\; \textbf{u}_t = 0&\;\text{ in }\; \mathbb{R}_{\geq 0}\times \Omega,\\
        p_0 = \bar{p}&\;\; \text{ on }\Omega,
    \end{split}
\end{equation}
where $\textbf{u}_t = (u^1_t,u^2_t):\mathbb{R}^2\to\mathbb{R}^2$ and $p_t:\mathbb{R}^2\to\mathbb{R}$ are time dependent functions representing velocity and pressure respectively. In the above system $\nabla^{\perp}p_t = (\partial_2 p_t,-\partial_1 p_t)$, represents a vector perpendicular to $\nabla p_t$ with the same norm. Due to various limitations while trying to attack the above problem directly, one is led to consider the following dual problem
\begin{equation}\label{dualSGE}
    \begin{split}
        \partial_t\rho_t + \text{div}(( y-\nabla P^*_t)^{\perp} \rho_t)=0&\;\text{ in }\; \mathbb{R}_{\geq 0}\times \mathbb{R}^2,\\
        \rho_t = (\nabla P_t)\#dx&\;\text{ in }\; \mathbb{R}_{\geq 0}\times \mathbb{R}^2,\\
        P_0=\bar{p}+\frac{|x|^2}{2}&\;\text{ on }\;\mathbb{R}^2.
    \end{split}
\end{equation}
Here $P_t = p_t+|x|^2/2$ is assumed to be convex and $P_t^*$ is the Legendre Transform of $P_t$ for each fixed time $t$. For a vector $v = (v_1,v_2)$, $v^{\perp} = (-v_2,v_1)$. One recovers the original velocity as,
\begin{equation}\label{vel_recover}
    \textbf{u}_t(x) = [\partial_t\nabla P^*_t](\nabla P_t(x))+D^2P^*_t(\nabla P_t(x))(\nabla P_t(x)-x)^{\perp}
\end{equation}
The long time existence for the dual problem had been established by Benamou-Brenier in \cite{BB1}. However while trying to recover $\textbf{u}_t$, the only apriori information one has about the term $D^2P^*_t$ is measure theoretic. This was resolved by Figalli-De Philippis in \cite{DF},  Figalli-De Philippis-Savin in \cite{DFS} and Schmidt in \cite{Schimdt} by deriving $W^{2,1+\varepsilon}$ estimates for Monge-Ampere equations with general right hand sides. The first existence result about distributional solutions to the semigeostrophic equation were presented in Ambrosio-Colombo-De Philippis-Figalli \cite{ACDF1}. In a followup paper \cite{ACDF2} by the same authors, they established similar results in 3 dimensions. \\

We always assume from now on that the measure $\rho_t$ is absolutely continuous w.r.t. Lebesgue measure on $\mathbb{R}^2$. We first note that if $P^*$ is a solution to the dual problem then by modding out by a constant rotation, the resulting profile satisfies a similar equation.

\begin{lemma}
    Let $(P^*,\rho)$ satisfy the dual Semigeostrophic equations \eqref{dualSGE}. Then $u(y,t) = P^*(R(t)y,t)$ satisfies
    \begin{equation}\label{new_sge}
    \begin{aligned}
        \partial_t \tilde{\rho} -\text{div}((\nabla u)^{\perp}\tilde{\rho})& = 0,\\
        \det(D^2u) = \tilde{\rho},
    \end{aligned}
    \end{equation}
for $\tilde{\rho}(y,t) = \rho (R(t)y,t)$. Here $R(t)\in SO(2)$ is the constant speed rotation $R(t) = e^{it}$.
\end{lemma}
\begin{proof}
    Note that
    $$\nabla u = R(-t) \nabla P^*,\qquad \det(D^2u) = \det(D^2P^*) = \tilde{\rho},$$
    for $\tilde{\rho}(y,t) = \rho(R(t)y,t)$. Therefore,
    $$\partial_t\rho = \partial_t\tilde{\rho} -\nabla \tilde{\rho}\cdot R(\frac{\pi}{2})R(-t)y.$$
    Calculating the divergence term in \eqref{dualSGE} we get,
    \begin{equation}
    \begin{aligned}
        \text{div}_y[(y-\nabla P^*)^{\perp}\rho]& = \text{div}_y[(y-R(t)\nabla u_t (R(-t)y))^{\perp}\tilde{\rho}(R(-t)y)]\\
        &=-\text{div}_y[(R(t)\nabla u(R(-t)y))^{\perp}\tilde{\rho}(R(-t)y)]+(R(-t))^{T}\nabla \tilde{\rho} \cdot y^{\perp} \\
        &= -\text{div}_y[(R(t)\nabla u(R(-t)y))^{\perp}\tilde{\rho}(R(-t)y)]+\nabla \tilde\rho \cdot R(-t)R(\frac{\pi}{2})y
    \end{aligned}
    \end{equation}
    Since rotations commute with each other. Therefore substituting into \eqref{dualSGE} we get,
    \begin{equation}
        \partial_t\tilde{\rho} -\text{div}_y((R(t)\nabla u(R(-t)y))^{\perp}\tilde{\rho}) = 0.
    \end{equation}
    Since for a vector field $F$, $\text{div}_{Ry}(F(Ry)) = \text{div}_{y}(R^{T}F(Ry))$ and rotations commute we get \eqref{new_sge}.
\end{proof}

Therefore if $P^*$ was a constantly rotating solution we get,
 \begin{equation}\label{rot_eqn}
    \begin{split}
        \text{div}((\nabla u)^{\perp} \tilde{\rho}) = 0\\
        \det(D^2u)=\tilde{\rho}
    \end{split}
 \end{equation}

 We call such solutions rotating wave solutions. All possible $(u,\tilde{\rho})$ that satisfy the above equations are of the form $\tilde{\rho} = \tilde{\rho}(u)$ at nondegenerate points of $u$ locally. We can in fact reconstruct $\rho$ globally,

\begin{lemma}\label{dual_SGE_thm}
    Let $P^*_t$ be a rotating wave solution to the dual SGE \eqref{dualSGE} with $\rho_t\not\equiv 0$. Then after a constant background rotation $R(t) = e^{it}\in SO(2)$, the solution to \eqref{rot_eqn} corresponding to $P^*$ is of the form $(u,f(u))$ satisfying,
    $$\det(D^2u) = f(u).$$
\end{lemma}

The assumption $\rho_t\not\equiv 0$ ensures $P^*$ cannot have a line in any sublevel set $K_s = \{P^*_t < s\}$ (see proof of Lemma \ref{no_line}). This forces the boundary of every level set $\partial K_s$ to be connected. Therefore away from the set $K_0 = \{P_t^* = \min\{P_t^*\}\}$ (if $\min\{P_t^*\}$ is defined) we always have $\nabla P_t^* \neq 0$. Moreover on the set $K_0$ we can define $\rho \equiv 0$.\\

 We have the following local boundedness of velocity for constantly rotating solutions to the SGE equation.

 \begin{theorem}\label{Thm_SGE1}
     Let $p_t$ be a rotating wave solution to the SGE equations \eqref{SGE} such that the corresponding $\rho_t = (\nabla P_t)\#dx$ is locally bounded between positive constants. Then the velocity $\textbf{u}_t$ is locally bounded.
 \end{theorem}
\begin{proof}
    The hypothesis of the above theorem ensures $P^*_t$ will be strictly convex in the corresponding dual coordinate region $\nabla P_t(\Omega)$. We can therefore locally always find a section of $P^*_t$ to apply Lemma \ref{dual_SGE_thm} and conclude using Theorem \ref{Main} and the formula \eqref{vel_recover}.
\end{proof}

 \begin{rmk}
     Given $(u,\rho(u))$ that satisfy conditions in Remark \ref{rmk:gen}, one can always construct constantly rotating solutions $P^*(x,t) = u(R(-t)x)$ to the dual SGE \eqref{dualSGE}. The velocity corresponding to this pair is always locally bounded and gives a general class of rotating wave solutions to the SGEs. We also note the previous work of McCann-Oberman in \cite{MO} where they construct solutions to SGE in a uniformly rotating elliptic tank with variable speed.
 \end{rmk}

 We also obtain the following global rigidity statement as well following from Theorem \ref{Liouville2},
 \begin{theorem}\label{SGE_thm2}
     Let $p_t$ be a rotating wave solution to the SGE equations \eqref{SGE} on $\mathbb{R}^2$ such that $ \nabla P_t(\mathbb{R}^2) = \mathbb{R}^2$. Then modulo a constant rotation, $\tilde{\rho}_t = f(P^*_t)$, for some nonnegative $f$. Further assume $f$ satisfies the conditions of Theorem \ref{Liouville2}. Then modulo a constant rotation, $P_t$ is affine equivalent to a radial function. 
 \end{theorem}

 Therefore if $P^*_t$ is a rotating wave solution, assume $P^*_t = u(R(-t)y)= u(z)$. Then using \eqref{vel_recover} we get,
 \begin{equation}\label{velocity_equation}
     \textbf{u}_t(R(t)\nabla u(z)) = R(t)(I-D^2u(z))R(\frac{\pi}{2})\nabla u(z).
 \end{equation}
 Therefore we have the relation,
 $$\textbf{u}_t(R(t)x) = R(t)\textbf{u}_0(x).$$
 Since $\textbf{u}_0$ is divergence free, there exists a $\psi$ such that $\textbf{u}_0 = R(\pi/2)\nabla \psi$. Therefore we have,
 $$\textbf{u}_t(x) = R(\frac{\pi}{2})\nabla_x [\psi(R(-t)x)]$$
Now consider the trajectory of a particle $X(t)$. We do a constant (reverse) rotation of coordinates $Z(t) = R(-t)X(t)$ to get,
\begin{equation}\label{trajectory}
    \frac{d}{dt}Z = -R(\frac{\pi}{2})Z +\textbf{u}_0(Z) = \left(\psi(Z)-\frac{1}{2}|Z|^2\right)^{\perp}
\end{equation}
\\
Therefore in rotating coordinates the particle trajectory lies along the level curves of,
$$H(Z) = \psi(Z)-\frac{1}{2}|Z|^2.$$
For a global rotating wave solution to the SGE \eqref{SGE} we have from Theorem \ref{SGE_thm2} that $u$ has to be affine equivalent to a radial function, i.e.
$$u(z) = \phi(q(z)), \quad q(z) = \frac{1}{2}z^TAz, \quad A \in \text{Sym}_{2 \times 2}, \quad A > 0.$$
It then follows from Theorem \ref{SGE_thm2} and equation \eqref{velocity_equation} that $\nabla H(x)$ is parallel to $Bx,\, B = R\left(-\frac{\pi}{2}\right)AR\left(\frac{\pi}{2}\right)$, i.e. that the level sets of $H$ are ellipses.

\begin{figure}[h!]
    \centering
    \includegraphics[width=0.7\linewidth]{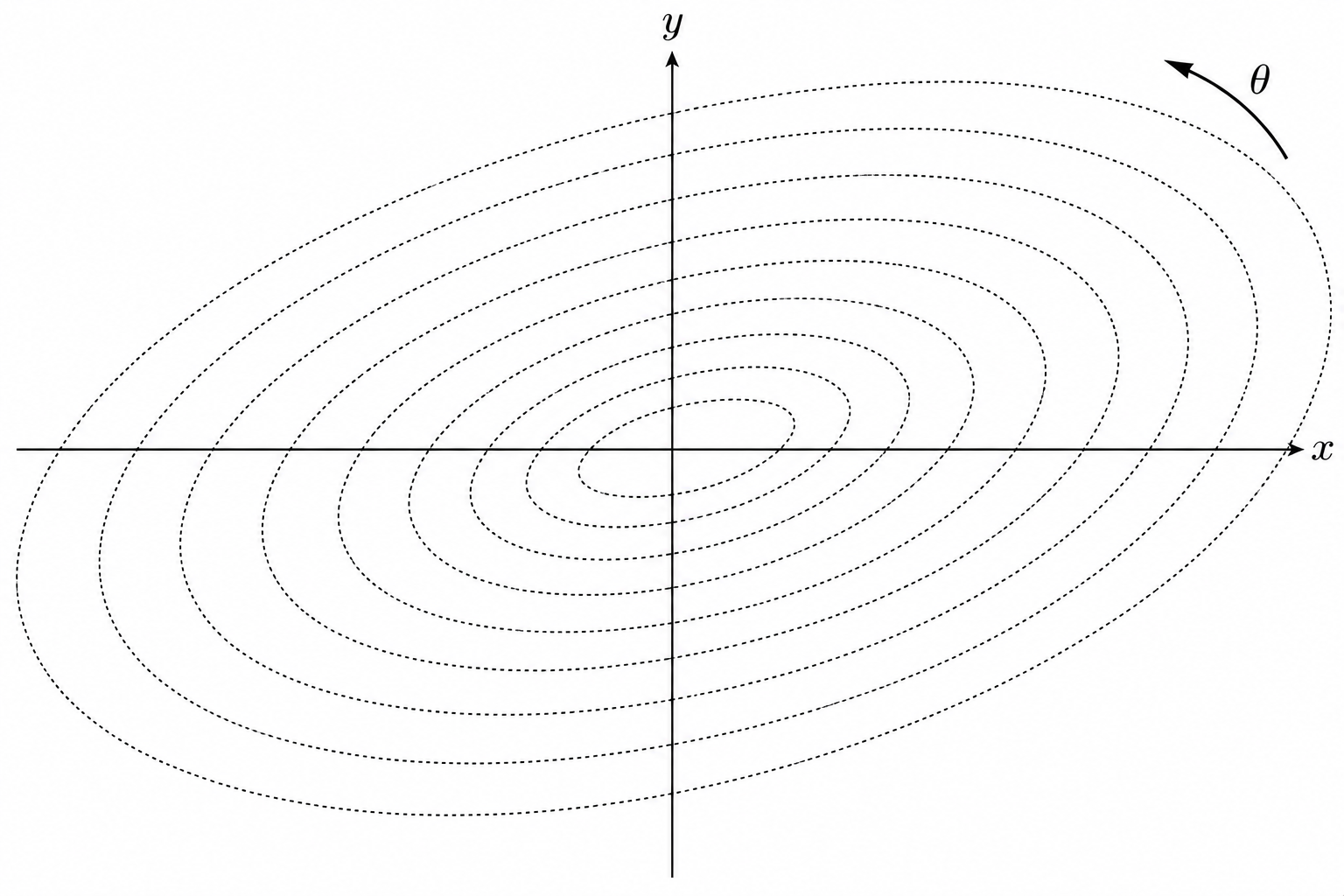}
    \caption{The trajectories in rotating coordinates.}
\end{figure}

\begin{figure}
    \centering
    \includegraphics[width=\linewidth]{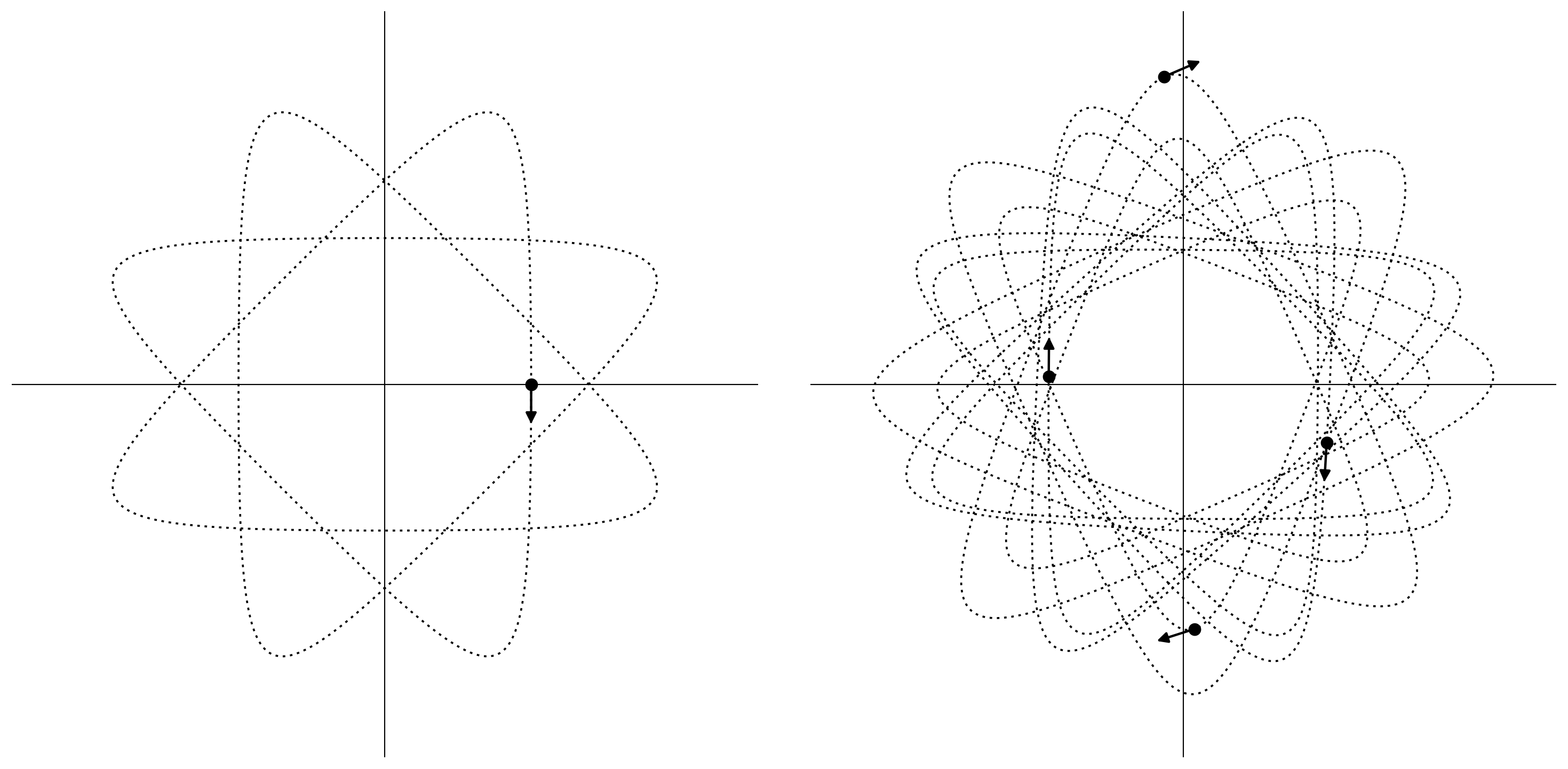}
    \caption{Trajectories of a single particle (left) and multiple particles starting at different points (right).}
    \label{fig:trajectory}
\end{figure}

The trajectory of each particle is governed by \eqref{trajectory}. In physical coordinates the path of a particle therefore looks like Figure \ref{fig:trajectory}.
\pagebreak
\section{Further Directions}\label{OpenQuestions}

Here we list a few natural questions inspired by this work.

\vspace{3mm}

\begin{enumerate}
\item Can Theorem \ref{Main} be extended to higher dimensions? If so, it may lead to a complete classification of global solutions to the obstacle problem \eqref{FB} (with $0 \leq p < n$) in all dimensions.

\vspace{5mm}

\item What can be said about the stability of rotating-wave solutions to the SGEs?

\vspace{5mm}

\item What can we say about the long-time behavior of solutions to the SGEs as a consequence of Theorem \ref{Main}?

\end{enumerate}

\end{document}